\documentclass[12pt]{article}

\usepackage{amssymb}
\usepackage{amsmath}
\usepackage{cite}
\usepackage[arrow, matrix, curve]{xy}
\usepackage{bbm}

\usepackage[dvips]{graphicx}
\usepackage{xcolor}

\begin{document}

\renewcommand{\citeleft}{{\rm [}}
\renewcommand{\citeright}{{\rm ]}}
\renewcommand{\citepunct}{{\rm,\ }}
\renewcommand{\citemid}{{\rm,\ }}

\newcounter{abschnitt}
\newtheorem{satz}{Theorem}
\newtheorem{coro}[satz]{Corollary}
\newtheorem{theorem}{Theorem}[abschnitt]
\newtheorem{koro}[theorem]{Corollary}
\newtheorem{prop}[theorem]{Proposition}
\newtheorem{lem}[theorem]{Lemma}
\newtheorem{expls}[theorem]{Examples}

\newcommand{\mres}{\mathbin{\vrule height 1.6ex depth 0pt width 0.11ex \vrule height 0.11ex depth 0pt width 1ex}}

\renewenvironment{quote}{\list{}{\leftmargin=0.62in\rightmargin=0.62in}\item[]}{\endlist}

\newcounter{saveeqn}
\newcommand{\alpheqn}{\setcounter{saveeqn}{\value{abschnitt}}
\renewcommand{\theequation}{\mbox{\arabic{saveeqn}.\arabic{equation}}}}
\newcommand{\reseteqn}{\setcounter{equation}{0}
\renewcommand{\theequation}{\arabic{equation}}}

\hyphenpenalty=9000

\sloppy

\phantom{a}

\vspace{-1.7cm}

\begin{center}
\begin{Large} {\bf Iterations of Minkowski Valuations} \\[0.4cm] \end{Large}

\begin{large} Oscar Ortega-Moreno \end{large}
\end{center}

\vspace{-0.9cm}

\begin{quote}
\footnotesize{ \vskip 1cm \noindent {\bf Abstract.}
It is shown that for any sufficiently regular even Minkowski valuation $\Phi$ which is homogeneous and intertwines rigid motions, and for any convex body $K$ in a smooth neighborhood of the unit ball, there exists a sequence of positive numbers $(\gamma_m)_{m=1}^\infty$ such that $(\gamma_m\Phi^m K)_{m=1}^\infty$ converges to the unit ball with respect to the Hausdorff metric.
}
\end{quote}

\vspace{0.4cm}

\centerline{\large{\bf{ \setcounter{abschnitt}{1}
\arabic{abschnitt}. Introduction}}}

\alpheqn

\vspace{0.5cm}

In the theory of convex bodies, a central role is played by the so-called projection body operator. To states its precise definition, let us recall that a convex body $K$ (that is, a compact, convex set with non-empty interior) in $\mathbb{R}^n$ is uniquely determined by its support function $h_K(u) = \max \{u \cdot x: x \in K\}$, $u \in  \mathbb{S}^{n-1}$. For a given convex body $K$, the support function of the projection body of $K$, denoted by $\Pi K$, is defined by $h_{\Pi K}(x) = \textrm{vol}_{n-1}\left(K|x^\perp\right)$ for all $x\in \mathbb{S}^{n-1}$, where $K|x^\perp$ is the projection of $K$ onto the hyperplane with normal unit vector $x$. The projection body was introduced by Minkowski himself and was later discovered to be of great importance in a number of different areas (see, e.g.,\textbf{\cite{gardner2ed}}, \textbf{\cite{KoldobskyBook}}, and \textbf{\cite{schneider93}}).

Petty's conjectured volume inequality for projection bodies  \textbf{\cite{Petty1972}} remains a major open problem in convex and integral geometry. It states that the quotient $V_n(\Pi K)/V_n(K)^{n-1}$ is minimized when $K$ is an ellipsoid. A positive answer to this conjecture would not only lead to a new isoperimetric inequality for convex sets, but also to an inequality stronger than a number of old and new isoperimetric inequalities, including important inequalities such as the polar projection inequality of Petty and the affine isoperimetric inequality.

The class reduction technique, first introduced by Schneider in \textbf{\cite{Schneider1987}}, reduces Petty's conjectured inequality to the study of a fixed point problem. More precisely, this approach consists in finding solutions to the equation $\Pi^2 K = \alpha K$ where $\alpha > 0$. This condition characterizes possible minimizers of the quotients $V_{n}(\Pi K)/V_{n}(K)^{n-1}$ and narrows the search for such minimizers down to the class of zonoids, which is the range of the projection body operator. 

As Lutwak observed in \textbf{\cite{Lutwak1990b,Lutwak1990a}}, the class reduction technique can be generalized to projection bodies of different degrees. For $i\in\{1,\dots, n-1\}$ and a given convex body $K$, the support function of the projection body of $K$ of order $i$, denoted by $\Pi_i K$, is defined by $h_{\Pi_i K}(x) = V_i\left(K|x^\perp\right)$ for all $x\in \mathbb{S}^{n-1}$, where $V_i$ denotes the $i$-th intrinsic volume. For $i\in \{1,\dots,n-1\}$, Lutwak \textbf{\cite{Lutwak1990b}} conjectured that the quotient $V_{i+1}(\Pi_i K)/V_{i+1}(K)^i$ is minimized among all convex bodies, when $K$ is a Euclidean ball. In fact, he showed that these chain of inequalities would follow from Petty's conjectured inequality. Furthermore, he extended the class reduction technique to characterize minimizers for these quotients by proving that a minimizer must also satisfy the fixed point equation $\Pi_i^2 K = \alpha K$ for some $\alpha>0$. In conjunction with a result by Schneider \textbf{\cite{Schneider1977}}, who showed that the only solutions to  $\Pi_1^2 K = \alpha K$ are balls, Lutwak \textbf{\cite{Lutwak1990b}} confirmed his conjecture for the case $i=1$ by showing that the quotient $V_{2}(\Pi_1 K)/V_{2}(K)$ is bounded from below.

While the solution to the equation $\Pi_1^2 K = \alpha K$ are only balls, it comes as a surprised that the projection body operator $\Pi = \Pi_{n-1}$ admits a variety of smooth and non-smooth solutions, with balls and cubes being the most basic examples. A complete characterization of polytopal solutions to the fixed point problem $\Pi^2 K = \alpha K$ was provided by Weil \textbf{\cite{Weil1971}}, who showed that they are the orthogonal sums of symmetric polygons and segments. Surprisingly, there is not much known about the smooth solutions. Recently, Ivaki \textbf{\cite{Ivaki2017,Ivaki2018}} made a major step towards this direction by showing that locally around the unit ball, the only smooth solution to the problem are ellipsoids. This also follows from a stronger result by Saroglou and Zvavitch~\textbf{\cite{SaroglouZvavitch}}, who considered limits of the sequence of iterations of the projection body operator with respect to the Banach-Mazur distance. With their result they manage to confirm Petty's conjectured inequality locally around the unit ball for smooth bodies. The main goal of this article is to extend the results of Saroglou and Zvavich~\textbf{\cite{SaroglouZvavitch}} on the convergence of the iterations of the projection body operator to a larger class of Minkowski valuations. A \emph{Minkowski valuation} on the space $\mathcal{K}^n$ of convex bodies in $\mathbb{R}^n$ endowed with the Hausdorff metric is a map $\Phi: \mathcal{K}^n \to \mathcal{K}^n$ such that
\begin{equation*}
\Phi K +\Phi L = \Phi(K\cup L) + \Phi(K\cap L)
\end{equation*}
whenever $K\cup L$ is convex and addition is the usual Minkowski addition.

The most celebrated theorem in valuation theory is undoubtedly Hadwiger's classic theorem \textbf{\cite{hadwiger51}} on the characterization of continuous, $\mathrm{SO}(n)$ invariant, real valued valuations as linear combinations of intrinsic volumes. It provided the kick-start needed to the systematic study of valuations and laid the groundwork for a number of venues of research in convex, discrete, and integral geometry (see, e.g., \textbf{\cite{Alesker1999, Alesker2001, Aleskerbook2018, BernigFu2011, HabParap2014, Klain:Rota, LudwigReitz2010}})

The investigation of Minkowski valuations is of a rather recent vintage. It started in 1974 with a paper by Schneider~\textbf{\cite{Schneider1974c}} on Minkowski endomorphisms: Minkowski additive valuations which commute with rigid motions. Nonetheless, it was the groundbreaking work by Ludwig  \textbf{\cite{ludwig02, Ludwig:Minkowski}} that prompted a wave of further developments. In her paper  \textbf{\cite{ludwig02}}, Ludwig solved a question posed by Lutwak on the characterization of the projection body operator as the unique continuous, translation invariant and  affinely contravariant Minkowski valuation. 

Following Ludwig's steps, many authors have confirmed a basic principle when it comes to characterizing  Minkowski valuations compatible with volume preserving linear maps: they often form a convex cone generated by finitely many transformations (see, e.g., \textbf{\cite{abardber11, boeroezkyludwig2019, colesantietal2017, haberl11, ludwig2010, SchuWan12, Wannerer2011}}). As opposed to affine compatiblitity, the less restrictive condition of being merely rotation equivariant means that the characterization of Minkowski valuation with such property will encompass a larger class of transformations: a cone with infinitely many generators. This feature makes the problem of finding an analogue of Hadwiger's theorem for Minkowski valuations a challenging one. Several authors have made substantial progress in this direction (see \textbf{\cite{kiderlen05, Schu06a, Schu10, SchuWan13, SchuWan16}}) by showing that Minkowski valuations admit a certain \emph{spherical convolution} representation (see Section 2). More recently, Dorrek \textbf{\cite{Dorrek2017b}} established the following remarkable  theorem under the additional assumption of homogeneity. Throughout, a map $\Phi: \mathcal{K}^n \rightarrow \mathcal{K}^n$ is said to have degree $i$ if $\Phi(\lambda K) = \lambda^i \Phi(K)$ for all $K \in \mathcal{K}$ and $\lambda \geq 0$. (By a result of McMullen, any translation invariant continuous valuation that is also homogeneous must be of integer degree $i \in \{0, \ldots, n\}$.)

\begin{theorem} \label{thmdorrek} \textnormal{(\!\!\!\textbf{\cite{Dorrek2017b}})}
If $\Phi_i: \mathcal{K}^n \rightarrow \mathcal{K}^n$ is a continuous translation invariant Minkowski valuation of degree $1 \leq i \leq n - 1$ which commutes with $\mathrm{SO}(n)$, then there exists a unique
 $\mathrm{SO}(n - 1)$ invariant $f \in L^1(\mathbb{S}^{n-1})$ with center of mass at the origin such that for every $K \in \mathcal{K}^n$,
\begin{equation} \label{convrep}
h_{\Phi_i K} = S_i(K,\cdot) \ast f.
\end{equation}
\end{theorem}

\vspace{0.2cm}

Here $S_i(K,\cdot)$, $1 \leq i \leq n - 1$, is the  \emph{area measure of order~$i$} on $\mathbb{S}^{n-1}$ associated to $K$ (see Section 2).
We call the function $f$ in (\ref{convrep}) the \emph{generating function} of $\Phi_i$. A complete classification of all such generating functions is still an important open problem in valuation theory. However, it is known that the support function of an arbitrary convex body of revolution $L \in \mathcal{K}^n$ is the generating function of a Minkowski valuation of degree $i$. In this case, we say that $\Phi_i$ is generated by $L$. If, in addition, for some $m\geq 2$, the boundary of $L$ is a $C^m$ submanifold of $\mathbb{R}^n$ with everywhere positive Gaussian curvature, we call $\Phi_i$ a \emph{$C^m_+$~regular} Minkowski valuation.
If $\Phi_i K = \{o\}$ for all $K \in \mathcal{K}^n$, we call $\Phi_i$ \emph{trivial}. 

In recent years, several classic inequalities involving projection bodies of arbitrary degree have been shown to hold for large (if not all) subclasses of Minkowski valuations intertwining rigid motions (see, e.g., \textbf{\cite{ABS2011, BPSW2014, HaberlSchu2019, HofstaetterSchu2021, papschu12, Schu06a}}). Some of these results are indeed a consequence of already known inequalities for the projection bodies, which turn out to be the limiting cases of such families of inequalities. In a recent article\textbf{\cite{OrtegaSchuster2021}}, the approach suggested by the author together with Schuster is the opposite one: showing that the inequality holds for a large and well behaved family of Minkowski valuations with which one can approximate the projection body operator.

\vspace{0.25cm}

\noindent {\bf Conjecture 1}\,\,(\!\!\!\textbf{\cite{OrtegaSchuster2021}}) Let \emph{$\Phi_i: \mathcal{K}^n \rightarrow \mathcal{K}^n$ be a continuous translation invariant Minkowski valuation of degree $1 \leq i \leq n - 1$ which commutes with $\mathrm{SO}(n)$. Among convex bodies $K \subseteq \mathbb{R}^n$ of non-empty interior, the volume ratio $V_{i+1}(\Phi_i K)/V_{i+1}(K)^{i}$ is minimized when $K$ is a Euclidean ball.}

\vspace{0.25cm}

As pointed out in \textbf{\cite{OrtegaSchuster2021}}, this conjecture provides a new approach to Petty's conjecture. First of all, notice that the projection body map can by approximated by $C^\infty_+$ regular Minkowski valuations. Thus, if the conjecture is proven for a sufficiently large class of these smooth valuations, we can obtain Petty's conjectured inequality (up to equality cases) by taking a limit. There are a couple of reasons why this smoothing approach seems advantageous. On one hand, if $\Phi_i$ is a smooth valuation, then for any convex body $K$, its image $\Phi K$ belongs to the same class of smoothness. This eliminates all possible polytopal solutions of the fixed point problem and thus reduces the set of possible minimizers of the quotient. On the other hand, the additional regularity of the operators allows us to use analytic machinery to study the problem (indeed, as we will see later, this is already reflected in the local solution around the unit ball).

A further generalization of the Schneider and Lutwak class reduction technique to Minkowski valuations was obtain in \textbf{\cite{OrtegaSchuster2021}}:

\vspace{0.3cm}

\noindent {\bf Proposition 1} (\!\!\textbf{\cite{OrtegaSchuster2021}}) \emph{Let $1 \leq i \leq n - 1$ and $\Phi_i: \mathcal{K}^n \rightarrow \mathcal{K}^n$ be a non-trivial monotone and translation invariant Minkowski valuation of degree $i$ which commutes with $\mathrm{SO}(n)$. If $K \in \mathcal{K}^n$ has dimension at least $i + 1$, then
\begin{equation} \label{prop1inequ}
\frac{V_{i+1}(\Phi_i K)}{V_{i+1}(K)^i} \geq \frac{V_{i+1}(\Phi_i^2 K)}{V_{i+1}(\Phi_iK)^i}
\end{equation}
with equality if and only if $\Phi_i^2 K$ and $K$ are homothetic. Moreover, if $\Phi_i$ is $C_2^+$~regular and balls are the only solutions to the fixed-point problem $\Phi_i^2 K = \alpha K$ for some $\alpha > 0$, then $V_{i+1}(\Phi_i K)/V_{i+1}(K)^i$ is minimized precisely by Euclidean balls.}

\vspace{0.3cm}

The problem to determine the convergence of iterations of Minkowski valuations is itself an intriguing one. Nevertheless, its importance becomes unquestionable once one realizes its connection to Conjecture~1 via the class reduction technique (Proposition~1). To explain this in more detail, let us first state the main result of this article.

\vspace{0.2cm}

\begin{satz} \label{mainthm} 
Let $2\leq i \leq  n-1$ and $\Phi_i: \mathcal{K}^n \to \mathcal{K}^n$ be a $C^\infty_+$ regular translation invariant even Minkowski valuation of degree $i$ which commutes with $\mathrm{SO}(n)$. Then, there exists $ \varepsilon >0$ such that if $K \in \mathcal{K}^n$ has support function of class $C^2$ and satisfies
$\|h_{\gamma K}-h_{\mathbb{B}^n}\|_{C^2(\mathbb{S}^{n-1})} < \varepsilon$
for some $\gamma>0$, then there exists a sequence of positive numbers $(\gamma_m)_{m=1}^\infty$ such that 
\begin{equation*}
d_{H}(\gamma_m\Phi_i^m K,\mathbb{B}^{n})\to 0 \quad \text{as }\,\, m \to \infty.
\end{equation*}
\end{satz}
\vspace{0.2cm}
Here $d_H(K,L)$ denotes the Hausdorff distance of two given sets $K$ and $L$.
%and $s_i(K,\cdot)$ is the density of the surface %area measure $S_i(K,\cdot)$ of order $i$ for all %$i\in {1,\dots, n-1}$. 
As a simple consequence of Theorem 1 we obtain a local solution to the fixed point problem $\Phi^2 K = \alpha K$ for smooth Minkowski valuations. This provides a different proof to the main results of \textbf{\cite{OrtegaSchuster2021}}. However, in \textbf{\cite{OrtegaSchuster2021}} the fixed point problem in question is solved for a bigger class of Minkowski valuations that also includes valuations that are not necessarily smooth. 

\vspace{0.3cm}

\noindent {\bf Corollary 1}\emph{ Let $2\leq i \leq  n-1$ and $\Phi_i: \mathcal{K}^n \to \mathcal{K}^n$ be a $C^\infty_+$ regular translation invariant even Minkowski valuation of degree $i$ which commutes with $\mathrm{SO}(n)$. Then there exists $ \varepsilon >0$ such that if $K \in \mathcal{K}^n$ has $C^2$ support function and satisfies
$\|h_{\gamma K}-h_{\mathbb{B}^n}\|_{C^2(\mathbb{S}^{n-1})} < \varepsilon$
for some $\gamma>0$, and 
\begin{equation*}
\Phi_i^2 K = \alpha K,
\end{equation*}
for some $\alpha>0$, then $K$ is a ball.
}

\vspace{0.3cm}

Furthermore, Theorem 1 leads to a local solution to Conjecture 1 using the class reduction technique. Note that Corollary 1 by itself does not suffice to obtain a local solution of the conjecture, as it only characterizes possible minimizers locally: balls can still be a different kind of critical points for the quotient. The proof relies heavily on the local convergence of the iterations. 

\vspace{0.3cm}

\noindent {\bf Corollary 2} \emph{ 
Let $2\leq i \leq  n-1$ and $\Phi_i: \mathcal{K}^n \to \mathcal{K}^n$be a $C^\infty_+$ regular translation invariant even Minkowski valuation of degree $i$ which commutes with $\mathrm{SO}(n)$. Then there exists $ \varepsilon >0$ such that if $K \in \mathcal{K}^n$ has $C^2$ support function and satisfies
$\|h_{\gamma K}-h_{\mathbb{B}^n}\|_{C^2(\mathbb{S}^{n-1})} < \varepsilon$
for some $\gamma>0$, then 
\begin{equation*}
\frac{V_{i+1}(\Phi_i K)}{V_{i+1}(K)^i}\geq \frac{V_{i+1}(\Phi_i \mathbb{B}^n)}{V_{i+1}(\mathbb{B}^n)^i},
\end{equation*} 
with equality if and only if $K$ is a ball.
}

\vspace{0.3cm}

The reader might be wondering why the case $i = 1$ is excluded from the statement of the previous theorems. This is due to the fact that, in this case, the convergence holds globally on $\mathcal{K}^{n}$.

\begin{satz} \label{mainthm2} 
Let $\Phi_1: \mathcal{K}^n \to \mathcal{K}^n$ be a translation invariant monotone Minkowski valuation of degree $1$ which commutes with $SO(n)$ with generating function $g$ of class $C^2$. Then, for all $K \in \mathcal{K}^n$
\begin{equation*}
d_{H}\left(\gamma_m \Phi_1^m K,\mathbb{B}^{n}\right)\to 0 \quad \text{as }\,\, m \to \infty,
\end{equation*}
with $\gamma_m^{-1} = \frac{w(K)}{2}\left(\int_{\mathbb{S}^{n-1}}g(u)du\right)^m$, where $w(K)$ denotes the mean width of $K$.
\end{satz}

The paper is organized as follows. In Section 2, we recall some basic background material from convex geometry, harmonic analysis and approximation theory. In Section 3, we state and prove some auxiliary results that will be needed for the proof of our main results. Finally, in the last section we prove the main theorem of this article as well as its corollaries. 

\vspace{.8cm}

\centerline{\large{\bf{ \setcounter{abschnitt}{2}
\arabic{abschnitt}. Background material}}}

\reseteqn \alpheqn \setcounter{theorem}{0}

\vspace{0.4cm}

The aim of this section is to recall some terminology and notation as well as standard results from convex geometry and harmonic analysis. Excellent and comprehensive presentations of the material cited here are the monographs by Gardner \textbf{\cite{gardner2ed}}, Schneider \textbf{\cite{schneider93}}, and Groemer \textbf{\cite{Groemer1996}}, with the latter focusing specifically on harmonic analysis and its applications to convex geometry.\\

 \noindent {\bf Convex bodies and Support functions.} A convex body is a convex compact subset of $\mathbb{R}^n$  with non-empty interior. We denote by $\mathcal{K}^n$ the set of all convex bodies in $\mathbb{R}^n$. Each convex body $K \in \mathcal{K}^n$ is determined by its support function 
\[h_K(x) = \max \{x \cdot y: y \in K\},\quad x \in \mathbb{R}^n,\]
which is (positively) homogeneous of degree one and subadditive. Conversely, every function on $\mathbb{R}^n$ with these properties is the support function of a unique convex body in $\mathcal{K}^n$. In particular, a homogeneous function $h \in C^2(\mathbb{R}^n)$ of degree one is the support function of a convex body $K \in \mathcal{K}^n$ if and only if its Hessian $\nabla^2h(u)$ is positive semi-definite for all $u \in \mathbb{S}^{n-1}$. Note that the gradient of a $1$-homogeneous function is a $0$-homogeneous map and so $\nabla^2f(u)u = 0$ for all $u\in \mathbb{S}^{n-1}$. The Hessian of $f$ is therefore determined by the restriction of the map $\nabla^2 f(u)$ to $u^\perp$. We will always refer to this restriction as $D^2f(u)$ for all $u\in \mathbb{S}^{n-1}$.

The Minkowski sum of $K, L \in \mathcal{K}^n$ is defined by $$K + L = \{x + y: x \in K, y \in L\}.$$ It is easy to check that the support function of $K+L$ is given by $h_{K + L} = h_K + h_L.$ For every $\vartheta \in \mathrm{SO}(n)$ and $y \in \mathbb{R}^n$, we have
\begin{equation*} \label{sontranslsuppfct}
h_{\vartheta K}(x) = h_K(\vartheta^{-1}x) \qquad \mbox{and} \qquad h_{K + y} = h_K(x) + x \cdot y
\end{equation*}
for all $x \in \mathbb{R}^n$. Moreover, $K \subseteq L$ if and only if $h_{K} \leq h_{L}$, in particular, $h_{K} > 0$ if and only if $o \in \mathrm{int}\,K$.\\

{\noindent\bf Higher regularity.} A body $K \in \mathcal{K}^n$ is said to be of class $C^k_+$ if its boundary hypersurface $\partial K$ is a $C^k$ submanifold of $\mathbb{R}^n$ and the map $n_K:\partial K \to \mathbb{S}^{n-1}$ that maps a boundary point to its unique outer unit normal is a $C^k$ diffeomorphism. Equivalently, $K \in \mathcal{K}^n$ is of class $C^k_+$ if $h_K \in C^k(\mathbb{R}^n)$ and the restriction of the Hessian $\nabla^2h_K(u)$ to $u^{\perp}$ is positive definite for every $u \in \mathbb{S}^{n-1}$ i.e, $D^2h_K(u)$ is positive definite for every $u \in \mathbb{S}^{n-1}$. We mention here that for any smooth function $f:\mathbb{S}^{n-1}\to \mathbb{R}$, we denote by $\nabla^2 f(u)$ the Hessian of its 1 homogeneous extension and $D^2 f(u)$ the restriction of $\nabla^2 f(u)$ to $u^\perp$. Moreover, if $\nabla^2_\mathbb{S}f$ denotes the spherical Hessian, then $D^2 f = \nabla^2_\mathbb{S}f + f I$, where $I$ denotes the identity.\\

{\noindent\bf Mixed volumes and area measures.} One of the pillars of the Brunn–Minkowski theory of convex bodies  is the fact that the volume of a Minkowski linear combination $\lambda_1K_1 + \cdots + \lambda_mK_m$, where $K_1, \ldots, K_m \in \mathcal{K}^n$ and $\lambda_1, \ldots,
\lambda_m \geq 0$, can be expressed as a homogeneous polynomial of degree $n$,
\begin{equation*} \label{mixed}
V_n(\lambda_1K_1 + \cdots +\lambda_m K_m)=\sum \limits_{j_1,\ldots, j_n=1}^m V(K_{j_1},\ldots,K_{j_n})\lambda_{j_1}\cdots\lambda_{j_n},
\end{equation*}
where the coefficients $V(K_{j_1},\ldots,K_{j_n})$ are the {\it mixed volumes} of $K_{j_1}, \ldots, K_{j_n}$, which depend only on $K_{j_1}, \ldots, K_{j_n}$ and are symmetric in their
arguments. Moreover, mixed volumes are translation invariant, Minkowski additive, monotone w.r.t.\ set inclusion in each of their arguments, and $V(K_1,\ldots,K_n) > 0$ if and only if there are segments
$l_i \subseteq K_i$, $1 \leq i \leq n$, with linearly independent directions.

For $K, L \in \mathcal{K}^n$ and $0 \leq i \leq n$, let $V(K[i],L[n-i])$ denote the mixed volume with $i$ copies of $K$ and $n - i$ copies of $L$. The \emph{$i$th intrinsic volume} of $K$ is given by
\begin{equation*} \label{viwi}
V_i(K)=\frac{1}{\kappa_{n-i}}\binom{n}{i} V(K[i],\mathbb{B}^n[n-i]),
\end{equation*}
where $\kappa_{m}$ denotes the $m$-dimensional volume of $\mathbb{B}^m$.

Associated with an $(n-1)$-tuple of bodies $K_2,\ldots, K_n \in \mathcal{K}^n$ is a finite Borel measure $S(K_2,\ldots,K_n,\cdot)$ on
$\mathbb{S}^{n-1}$, the {\it mixed area measure}, such that for all $K_1 \in \mathcal{K}^n\!$,
\begin{equation} \label{defmixedarea}
V(K_1,\ldots,K_n)=\frac{1}{n}\int\displaylimits_{\mathbb{S}^{n-1}} h(K_1,u)\,dS(K_2,\ldots,K_n,u).
\end{equation}

For $K \in \mathcal{K}^n$ and $0 \leq i \leq n - 1$, the measures $S_i(K,\cdot) := S(K[i],\mathbb{B}^n[n-1-i],\cdot)$ are called the \emph{area measures} of order $i$ of $K$. The measure $S_{n-1}(K,\cdot)$ is
also known as the \emph{surface area measure} of $K$. If $K$ has non-empty interior, then, by a theorem of Aleksandrov--Fenchel--Jessen (see, e.g.,
\textbf{\cite[\textnormal{p.\ 449}]{schneider93}}), each of the measures $S_i(K,\cdot)$, $1 \leq i \leq n - 1$, determines $K$ up to translations.
The centroid of every area measure of a convex body is at the origin, that is, for every $K \in \mathcal{K}^n$ and all $i \in \{0, \ldots, n - 1\}$,
\begin{equation*}
\int\displaylimits_{\mathbb{S}^{n-1}} u\,dS_i(K,u) = o.
\end{equation*}
\emph{Minkowski's existence theorem} states that a non-negative Borel measure $\mu$ on $\mathbb{S}^{n-1}$ is the
surface area measure of some $K \in \mathcal{K}^n$ with non-empty interior if and only if $\mu$ is not concentrated on a great subsphere of $\mathbb{S}^{n-1}$ and has centroid at the origin.\\

{\noindent\bf Area densities and mixed discriminants.} If $K \in \mathcal{K}^n$ has a $C^2$ support function, then each measure $S_i(K,\cdot)$, $0 \leq i \leq n - 1$, is absolutely continuous w.r.t.\ spherical Lebesgue measure. To make this more precise, let us recall the notion of mixed discriminants.
If $A_1, \ldots, A_m$ are symmetric real $k \times k$ matrices and $\lambda_1,\dots,\lambda_m \geq 0$, then
\begin{equation}\label{polydet}
\det(\lambda_1A_1+\cdots+\lambda_m A_m) = \sum_{j_1,\dots,j_k = 1}^m \mathrm{D}(A_{j_1},\ldots,A_{j_k})\lambda_{j_1}\cdots\lambda_{j_k},
\end{equation}
where the coefficients $\mathrm{D}(A_{j_1},\ldots,A_{j_k})$ are the {\it mixed discriminants} of $A_{j_1}, \ldots, A_{j_k}$ which depend only on $A_{j_1}, \ldots, A_{j_k}$ and are symmetric and multilinear in their arguments. Clearly, $\mathrm{D}(A,\dots,A) = \det(A)$ for any symmetric $k \times k$ matrix $A$. Moreover, $\mathrm{D}(B A_1,\ldots, B A_k) = \det(B)\mathrm{D}(A_1,\ldots,A_k)$,
\begin{equation} \label{basicMD}
\mathrm{D}(A,B,\ldots,B) = \frac{1}{n-1}\mathrm{tr}(\mathrm{cof}(B)A)
\end{equation}
for any symmetric $k \times k$ matrix $B$, and if $A_1,\ldots,A_k$ are positive semi-definite, then $\mathrm{D}(A_1,\dots,A_k) \geq 0$. Finally, if $K_1,\ldots,K_{n-1} \in \mathcal{K}^n$ have support functions $h_1,\ldots,h_{n-1} \in C^2(\mathbb{R}^n)$, then the density of $S(K_1,\ldots,K_{n-1},\cdot)$ is given by
\begin{equation} \label{densofmixedarea}
s(K_1,\ldots,K_{n-1},u) = \mathrm{D}(D^2h_1(u),\ldots,D^2h_{n-1}(u)), \quad u \in \mathbb{S}^{n-1}.
\end{equation}
In particular, for $K \in \mathcal{K}^n$ with support function $h \in C^2(\mathbb{R}^n)$, we have
\begin{equation} \label{densofsurfarea}
s_{n-1}(K,u) = \det D^2h(u),\quad u \in \mathbb{S}^{n-1}.
\end{equation}
Motivated by (\ref{densofmixedarea}) and (\ref{densofsurfarea}), we frequently use in subsequent sections the notation $s(h_1,\ldots,h_{n-1},\cdot)$, $s_{n-1}(h,\cdot)$, $\ldots$ instead of
$s(K_1,\ldots,K_{n-1},\cdot)$, $s_{n-1}(K,u)$, etc.\\

{\noindent \bf Metrics on the space $\mathcal{K}^n$.} The most used notion of convergence on convex bodies is the one derived from the Hausdorff metric. Given any pair of non-empty sets $X,Y\subset \mathbb{R}^n$, their Hausdorff distance is defined by
$$d_H(X,Y) = \inf\{\varepsilon \geq 0\,;\ X \subseteq Y_\varepsilon \text{ and } Y \subseteq X_\varepsilon\},$$
where 
$$X_\varepsilon = \bigcup_{x \in X} \{z \in M\,;\ \|z- x\|_2 \leq \varepsilon\}.$$ It is well known that the Hausdorff metric of two bodies $K, L \in \mathcal{K}^n$ can be expressed as 
$$d_H(K,L) = \|h_K - h_L\|_{\infty},$$ 
where $\|\,\cdot\,\|_{\infty}$ denotes the maximum norm on $C(\mathbb{S}^{n-1})$.

Naturally, a further family of metrics is derived  by replacing the maximum norm by an $L^p$ norm. For $p\in [1,\infty)$ and  $K,L\in \mathcal{K}^n$, let
$$d_p(K,L) = \left(\int_{\mathbb{S}^{n-1}}|h_K(u) - h_L(u)|^p du\right)^\frac{1}{p}.$$
Clearly, $d_p$ is a metric on $\mathcal{K}^n$, called the $L_p$ metric. The $L_2$ metric is of particular interest in this context because of the well-developed Harmonic analysis on $\mathbb{S}^{n-1}$. In general, $L_p$ metrics do not necessarily induce the same notion of convergence; however, with the additional structure provided by the support functions, Vitale \textbf{\cite{Vitale1985}} deduced that all of the $d_p$ metrics, $p\in [1,\infty]$, generate the same topology on $\mathcal{K}^n$ and yield a complete metric spaces in which closed, bounded sets are compact (see Theorem 3 in \textbf{\cite{Vitale1985}}).

There are other distances on convex bodies that emerge from distances between measures. By the Aleksandrov--Fenchel--Jessen uniqueness theorem (see, e.g., \textbf{\cite[\textnormal{p.\ 449}]{schneider93}}), a convex body is uniquely determined, up to translation, by its area measure of order $i$ for all $i\in \{1,\dots,n-1\}$. Therefore, any metric in the space of positive measures $\mathcal{M}(\mathbb{S}^{n-1})$ gives rise to a metric on the set of convex bodies with centroid at the origin. Some examples of such metrics, that will play an important role in the sequel are the total variation and the L\'evy--Prokhorov metric. The total variation metric between to positive measures $\mu,\nu \in \mathcal{M}(\mathbb{S}^{n-1})$ is given by 
\[d_{\mathrm{TV}}(\mu,\nu)=\sup_{ A\in \mathcal{B}(\mathbb{S}^{n-1})}\left|\mu(A)-\nu(A)\right|\]
where $ \mathcal{B}(\mathbb{S}^{n-1})$ denotes the set of Borel sets on $\mathbb{S}^{n-1}$. The L\'evy--Prokhorov metric is a weaker metric on the space of positive measures on the sphere $\mathcal{M}_+(\mathbb{S}^{n-1})$. It is defined by 
\begin{equation*}
d_{\mathrm{LP}}(\mu,\nu) = \inf \left\{\varepsilon >0 : \mu(A)\leq \nu(A_\varepsilon) + \varepsilon \text{ and }  \nu(A)\leq \mu(A_\varepsilon) + \varepsilon, A \in \mathcal{B}(\mathbb{S}^{n-1}) \right\},
\end{equation*}
where $A_\varepsilon = \{u\in \mathbb{S}^{n-1}: d(A,u)<\varepsilon \}$. The metric $d_\mathrm{LP}$ has the property that $d_{\mathrm{LP}}(\mu_k, \mu) \to 0$ if and only if $(\mu_k)_{k=1}^\infty$ converges weakly to $\mu$.
It is easy to see that the convergence with respect to $d_{TV}$ implies convergence with respect to $d_{LP}$.
\\ 

{\noindent\bf Harmonic analysis on $\mathbb{S}^{n-1}$.} We turn now to the background material on spherical harmonics. To this end, let $\Delta_{\mathbb{S}}$ denote the spherical Laplacian on $\mathbb{S}^{n-1}$ and recall that it is
a second-order uniformly elliptic self-adjoint operator. We write $\mathcal{H}_k^n$ for the vector space of spherical harmonics of dimension $n$ and degree $k$ and denote its dimension by
\begin{equation} \label{nnk}
N(n,k) = \frac{n + 2k - 2}{n + k - 2} {n + k - 2 \choose n - 2} = \mathrm{O}(k^{n - 2}) \mbox{ as } k \rightarrow \infty.
\end{equation}
Spherical harmonics are (precisely) the eigenfunctions of $\Delta_{\mathbb{S}}$, more specific, for $Y_k \in \mathcal{H}_k^n$, we have
\begin{equation} \label{deltasmult}
\Delta_{\mathbb{S}} Y_k = -k(k + n - 2)\,Y_k.
\end{equation}

The spaces $\mathcal{H}_k^n$ are pairwise orthogonal subspaces of $L^2(\mathbb{S}^{n-1})$. Moreover, the Fourier series $f \sim \sum_{k=0}^{\infty} \pi_k f$
converges to $f$ in $L^2$ for every $f \in L^2(\mathbb{S}^{n-1})$, where $\pi_k: L^2(\mathbb{S}^{n-1}) \rightarrow \mathcal{H}_k^n$ denotes the orthogonal projection.
Letting $P_k^n \in C([-1,1])$ denote the \emph{Legendre polynomial} of dimension $n$ and degree $k$, we have
\begin{equation} \label{projleg}
(\pi_k f)(v) = \frac{N(n,k)}{\omega_n} \int\displaylimits_{\mathbb{S}^{n-1}}\!\! f(u)\, P_k^n(u\cdot v) \,du, \qquad v \in \mathbb{S}^{n-1},
\end{equation}
where $\omega_n$ denotes the surface area of $\mathbb{B}^n$ and integration is with respect to spherical Lebesgue measure. Since the orthogonal projection $\pi_k$ is self adjoint, it is consistent to extend it to the space $\mathcal{M}(\mathbb{S}^{n-1})$ of signed finite Borel measures by
\begin{equation*}
(\pi_k \mu)(v) = \frac{N(n,k)}{\omega_n}\int\displaylimits_{\mathbb{S}^{n-1}}\!\!P^n_k(u\cdot v)\, d\mu(u), \qquad v \in \mathbb{S}^{n-1}.
\end{equation*}
It can be shown easily that $\pi_k \mu \in \mathcal{H}_k^n$ for all $k \geq 0$ and that the formal Fourier series
$\mu \sim \sum_{k = 0}^\infty \pi_k \mu$ uniquely determines the measure $\mu$.

Throughout, we use $\bar{e} \in \mathbb{S}^{n-1}$ to denote a fixed but arbitrarily chosen pole of $\mathbb{S}^{n-1}$ and write $\mathrm{SO}(n-1)$ for the stabilizer
in $\mathrm{SO}(n)$ of $\bar{e}$. In the theory of spherical harmonics, a function or measure on $\mathbb{S}^{n-1}$ which is $\mathrm{SO}(n-1)$ invariant is often called \emph{zonal}.
Clearly, zonal functions depend only on the value of $u \cdot \bar{e}$. The subspace of zonal functions in $\mathcal{H}_k^n$ is $1$-dimensional for every $k \geq 0$ and spanned by $u \mapsto P_k^n(u \cdot \bar{e})$.
Since the spaces $\mathcal{H}_k^n$ are orthogonal, it is not difficult to show that any zonal measure $\mu \in \mathcal{M}(\mathbb{S}^{n-1})$
admits a series expansion of the form
\begin{equation} \label{expzonal}
\mu \sim \sum_{k=0}^{\infty} \frac{N(n,k)}{\omega_n}\, a_k^n[\mu]\,P_k^n(\,\,.\cdot \bar{e}),
\end{equation}
where
\begin{equation} \label{multleg}
a_k^n[\mu] = \omega_{n-1} \int\displaylimits_{-1}^1 P_k^n(t)\,(1-t^2)^{\frac{n-3}{2}}\,d\bar{\mu}(t).
\end{equation}
Here, we have used cylindrical coordinates $u = t\bar{e} + \sqrt{1 - t^2} v$ on $\mathbb{S}^{n-1}$ to identify the zonal measure $\mu$ with a measure $\bar{\mu}$ on $[-1,1]$. If $\mu$ is absolutely continuous with density $f$ w.r.t.\ spherical Lebesgue measure, we write $a_k^n[f\,]$ instead of $a_k^n[\mu]$.

For the explicit computation of integrals of the form (\ref{multleg}) the following \emph{Formula of Rodrigues} for the Legendre polynomials
is often useful:
\begin{equation}\label{RodriguesF}
P_k^n(t)= (-1)^k \frac{\Gamma\left(\frac{n - 1}{2}\right)}{2^k \Gamma\left( \frac{n - 1}{2} + k\right)} (1 - t^2)^{-\frac{n-3}{2}}\left(\frac{d}{dt}\right)^k(1 - t^2)^{\frac{n-3}{2}+k}.
\end{equation}

{\noindent\bf Multiplier and integral transforms.} The well known \emph{Funk--Hecke Theorem} states that if $f \in C([-1,1])$ and $\mathrm{T}_{f}: \mathcal{M}(\mathbb{S}^{n-1}) \rightarrow C(\mathbb{S}^{n-1})$ is defined by
\begin{equation} \label{Funkheckinttrafo}
(\mathrm{T}_{f}\mu)(u) = \int\displaylimits_{\mathbb{S}^{n-1}}\!\! f(u \cdot v)\,d\mu(v), \qquad u \in \mathbb{S}^{n-1},
\end{equation}
then the spherical harmonic expansion of $\mathrm{T}_{f}\mu$ is given by
\begin{equation} \label{funkhecke}
\mathrm{T}_{f}\mu \sim \sum_{k=0}^\infty a_k^n[f]\,\pi_k\mu,
\end{equation}
where the numbers $a_k^n[f]$ are given by (\ref{multleg}) and called the \emph{multipliers} of $\mathrm{T}_{f}$.\\

Integral transforms of the form (\ref{Funkheckinttrafo}) are closely related to the convolution between functions and measures on $\mathbb{S}^{n-1}$.
In order to recall its definition, first note that the convolution $\sigma \ast \tau$ of signed measures $\sigma, \tau$ on the compact Lie group $\mathrm{SO}(n)$ can be defined by
\[\int\displaylimits_{\mathrm{SO}(n)}\!\!\! f(\vartheta)\, d(\sigma \ast \tau)(\vartheta)=\int\displaylimits_{\mathrm{SO}(n)}\! \int\displaylimits_{\mathrm{SO}(n)}\!\!\! f(\eta \theta)\,d\sigma(\eta)\,d\tau(\theta), \qquad f \in C(\mathrm{SO}(n)).   \]
By identifying $\mathbb{S}^{n-1}$ with the homogeneous space $\mathrm{SO}(n)/\mathrm{SO}(n-1)$, one obtains a one-to-one correspondence of $C(\mathbb{S}^{n-1})$ and $\mathcal{M}(\mathbb{S}^{n-1})$ with right
$\mathrm{SO}(n-1)$ invariant functions and measures on $\mathrm{SO}(n)$, respectively. Using this correspondence, the convolution of measures on $\mathrm{SO}(n)$ induces a convolution product
on $\mathcal{M}(S^{n-1})$ (for more details see, e.g., \textbf{\cite{Schu06a}}). For this spherical convolution, zonal functions and measures play a particularly important role.
Let us therefore denote by $C(\mathbb{S}^{n-1},\bar{e})$ the set of continuous zonal functions on $\mathbb{S}^{n-1}$. Then, for $\mu \in \mathcal{M}(\mathbb{S}^{n-1})$, $f \in C(\mathbb{S}^{n-1},\bar{e})$, and $\eta \in \mathrm{SO}(n)$, we have
\begin{equation} \label{zonalconv}
(\mu \ast f)(\eta \bar{e}) =\int\displaylimits_{\mathbb{S}^{n-1}}\!\!f(\eta^{-1} u)\,d\mu(u).
\end{equation}

Note that, by (\ref{zonalconv}), we have $(\vartheta \mu) \ast f = \vartheta(\mu \ast f)$ for every $\vartheta \in \mathrm{SO}(n)$,
where $\vartheta\mu$ is the image measure of $\mu$ under $\vartheta \in \mathrm{SO}(n)$. Moreover, from the identification of a
zonal function $f$ on $\mathbb{S}^{n-1}$ with a function $\bar{f}$ on $[-1,1]$, (\ref{multleg}), and (\ref{zonalconv}), we obtain
\begin{equation} \label{multconvpkn}
a_k^n[f] = \int\displaylimits_{\mathbb{S}^{n-1}}\!\!f(u) P_k^n(\bar{e}\cdot u)\,du,
\end{equation}
and the Funk--Hecke Theorem implies that
\begin{equation} \label{convmulttransf}
\mu \ast f \sim \sum_{k=0}^{\infty} a_k^n[f]\,\pi_k\mu.
\end{equation}
Hence, convolution from the right induces a multiplier transformation. It is also easy to check from (\ref{zonalconv}) that the convolution of zonal functions and measures is Abelian and that for all $\mu, \tau \in \mathcal{M}(\mathbb{S}^{n-1})$ and every $f \in C(\mathbb{S}^{n-1},\bar{e})$,
\begin{equation} \label{convselfad}
\int\displaylimits_{\mathbb{S}^{n-1}}\!\! (\mu \ast f)(u)\,d\tau(u) = \int\displaylimits_{\mathbb{S}^{n-1}}\!\! (\tau \ast f)(u)\,d\mu(u).
\end{equation}

\vspace{0.1cm}

\begin{expls} \label{exps1} \end{expls}

\vspace{-0.2cm}

\begin{enumerate}
\item[(a)] The cosine transform $\mathrm{C}: \mathcal{M}(\mathbb{S}^{n-1}) \rightarrow C(\mathbb{S}^{n-1})$ is defined by
\begin{equation}
\mathrm{C}\mu(u) = \int\displaylimits_{\mathbb{S}^{n-1}}\!\!|u \cdot v|\,d\mu(v) = (\mu \ast |\bar{e}\cdot\,.\,|)(u), \qquad u \in \mathbb{S}^{n-1}.
\end{equation}
Using the Formula of Rodrigues, the multipliers $a_k^n[\mathrm{C}] := a_k^n[|\bar{e}\cdot\,.\,|]$ of the cosine transform can be easily computed to
\begin{equation} \label{multC}
a_k^n[\mathrm{C}] =  (-1)^\frac{n-2}{2}2\frac{1\cdot 3\cdots (k-3)}{(n+1)(n+3)\cdots (k+n-1)}
\end{equation}
for even $k$ and $a_k^n[\mathrm{C}]= 0$ for $k$ odd.

\item[(b)] Generalizing (a), we define for an arbitrary body of revolution $L \in \mathcal{K}^n$ the integral transform $\mathrm{T}_{\!L}: \mathcal{M}(\mathbb{S}^{n-1}) \rightarrow C(\mathbb{S}^{n-1})$ by
\[\mathrm{T}_L \sigma = \sigma \ast h_L.\]
We denote its multipliers by $a_k^n[L] := a_k^n[h_L]$.
\end{enumerate}

Finally, recall that the second order differential operator $\Box_n$, defined by
\[\Box_n h = h +  \frac{1}{n-1}\Delta_{\mathbb{S}}h \]
for $h \in C^2(\mathbb{S}^{n-1})$, relates the support function $h_K$ of a convex body $K \in \mathcal{K}^n$ with its first-order area measure $S_1(K,\cdot)$ by
\begin{equation} \label{boxhks1}
\Box_n h_K = S_1(K,\cdot),
\end{equation}
where (\ref{boxhks1}) has to be understood in a distributional sense if $h_K$ is not of class $C^2$. From the definition of $\Box_n$ and (\ref{deltasmult}), we see
that for $h \in C(\mathbb{S}^{n-1})$ and every $k \geq 0$,
\begin{equation} \label{boxnmult}
\pi_k \Box_n h = \frac{(1-k)(k+n-1)}{n-1} \pi_k h.
\end{equation}
In particular, $\Box_n$ acts as a multiplier transformation and since such operators clearly commute, we note for later quick reference that, by (\ref{convmulttransf}) and (\ref{boxnmult}), we have
\begin{equation} \label{boxconvcomm}
\mathrm{T}_{f}\,\Box_n = \Box_n \mathrm{T}_f.
\end{equation}

{\noindent \bf The $\mathcal{U}_\alpha$ classes.} For a bounded real-valued measurable function $f$ on $\mathbb{S}^{n-1}$ let $\|f\|_{\mathcal{U}_\alpha}$ denote the smallest constant $M$ such that $\|f\|_{\infty} \leq M $ and, such that for all positive integers $k$, there exists a polynomial $p_k$ of degree $k$ such that $\|f-p_k\|_{L^2}\leq M k^{-\alpha}$. We define the class of real-valued functions $\mathcal{U}_\alpha$ to be the subset of bounded functions $f$ such that $\|f\|_{\mathcal{U}_\alpha}<\infty$. Let $\Theta : [0,\infty) \to [0,1]$ be an infinitely smooth function on such that $\Theta = 1$ on $[0,1]$ and $\theta = 0$ on $[2,\infty)$. Define the multiplier operator $\mathrm{M}_k$ by 
\begin{equation}\label{Mn}
\mathrm{M}_j f = \sum_{k = 0}^\infty \Theta(k/j)\,\, \pi_k f. 
\end{equation}
Note that $\mathrm{M}_j$ is a polynomial of degree at most $2j$ and that $\mathrm{M}_j\, p = p$ for any polynomial $p$ of degree $j$. A well known fact from approximation theory tells us that the operators $\mathrm{M}_j$ are uniformly bounded in $L_p$, i.e. there exist a constant $C = C(\Theta)$ such that $\|\mathrm{M}_j\|_{L^p\to L^p}\leq C(\Theta)$. A complete proof of this statement can be found in Appendix A in \textbf{\cite{FNRZ2011}}. 

\begin{lem}[\!\!\cite{FNRZ2011}]\label{Uclassesmult0}
Let $\alpha \geq 0$. If $f\in\mathcal{U}_\alpha$, then 
\[\|f - \mathrm{M}_j f\|_{\mathcal{U}_\alpha} \leq C \|f\|_{\mathcal{U}_\alpha}j^{-\alpha},\]
for some $C =C_\alpha>0$.
\end{lem}

\begin{lem}[\!\!\cite{FNRZ2011}]\label{Uclassesmult}
Let $\alpha \geq 0$. If $f,g\in\mathcal{U}_\alpha$, then $fg\in \mathcal{U}_\alpha$ and 
\[\Vert fg \Vert_{\mathcal{U}_\alpha} \leq C \Vert f\Vert_{\mathcal{U}_\alpha}\Vert g\Vert_{\mathcal{U}_\alpha},\]
where $C = C_\alpha>0$. 
\end{lem}

\begin{lem}[\!\!\cite{FNRZ2011}]\label{tradeoff}
Let $\beta> \alpha$. For every $\delta > 0$ there exists a constant $C_{\delta,\alpha,\beta} >0$ such that
$\|f\|_{\mathcal{U}_\alpha} \leq C_{\delta,\alpha,\beta} \|f\|_{\infty} + \delta\|f\|_{\mathcal{U}_\beta}$.  
\end{lem}

\begin{lem}[\!\!\cite{FNRZ2011}]\label{fromL2toLinf}
There exist constants $\alpha, C_1, C_2 > 0$ with the following properties: if $\varphi:\mathbb{S}^{n-1}\to \mathbb{R}$ satisfies $\|\varphi\|_{2}<\varepsilon$ for some $\varepsilon \in (0,1)$, and $\|\varphi\|_{\mathcal{U}_\alpha}\leq C_1$, then $\|\varphi\|_{\infty}< C_2\varepsilon^\frac{4}{n+3}$. 
\end{lem}

\vspace{0.8cm}

\centerline{\large{\bf{ \setcounter{abschnitt}{3}
\arabic{abschnitt}. Auxiliary Results}}}
\reseteqn \alpheqn \setcounter{theorem}{0}

\vspace{0.6cm}
In this section, we gather some auxiliary results that will be needed in the proof of our main theorem. The next proposition is a generalization of Proposition~3.2 in \textbf{\cite{GAPS}} to area measures of all degrees.

\begin{prop}\label{DpH} Let $i\in \{1,\dots, n-1 \}$ and $(K_m)_{m=1}\infty$ be a sequence of convex bodies in $\mathcal{K}^{n-1}$ with centroid at the origin. Then, $d_{LP}(S_i(K_m,\cdot),S_i(K,\cdot))\to 0$ as $k\to \infty$ if and only if $d_H(K_m,K)\to 0$ as $k\to \infty$.
\end{prop}
\noindent {\it Proof.} Suppose that there exists a convex body with centroid at the origin $K$ such that $d_{LP}(S_i(K_m,\cdot),S_i(K,\cdot))\to 0$ as $m\to \infty$. Since convergence in the L\'evy-Prokhorov metric is equivalent to weak convergence of measures, 
\begin{equation*}
V_i(K_m) = \int\displaylimits_{\mathbb{S}^{n-1}}dS_{i}(K_m,v) \to \int\displaylimits_{\mathbb{S}^{n-1}}dS_{i}(K,v) = V_i(K) \quad \text{as } m\to \infty.
\end{equation*}
Hence, there exists a constant $C>0$ such that $V_i(K_m)\leq C$ for all $m$. By the isoperimetric inequality between consecutive intrinsic volumes, i.e. $V_{i+1}(K)^i V_i(\mathbb{B}^n)^i \leq V_i(K)^{i+1} V_{i+1}(\mathbb{B}^n)^i$, there exists a constant $C'>0$ such that $V_{i+1}(K_m)\leq C'$ for all $m$. Again by the weak convergence of the area measures of order $i$, for all $u \in \mathbb{S}^{n-1}$
\begin{equation}
 h_{\Pi_iK_m}(u) = \frac{1}{2}\int\displaylimits_{\mathbb{S}^{n-1}}|u\cdot v|dS_{i}(K_m,v) \to \frac{1}{2}\int\displaylimits_{\mathbb{S}^{n-1}}|u\cdot v|dS_{i}(K,v) = h_{\Pi_iK}(u)
\end{equation} 
as $m\to \infty$. Thus, the sequence $(h_{\Pi_i K_m})_{m=1}^\infty$ converges uniformly to $ h_{\Pi_iK}$ (see {\bf\cite[\textnormal{Theorem 1.8.12}]{schneider93}}). Since $\Pi_i K$ is a full dimensional convex body with centroid at the origin and the area measures of order $i$ have centroid at the origin, there exists a constant $C''>0$ such that 
\begin{equation}\label{LP1}
C''\leq h_{\Pi_i K_m}(u) =  \frac{1}{2}\int\displaylimits_{\mathbb{S}^{n-1}}|u\cdot v|dS_i(K_m,v) = \int\displaylimits_{\mathbb{S}^{n-1}}\max\{0,u\cdot v\}dS_i(K_m,v),
\end{equation}
for all $u\in \mathbb{S}^{n-1}$ and all $m$. 
On the other hand, for any $s>0$ and $u\in \mathbb{S}^{n-1}$ and $m$ with $[o,su]\subseteq K_m$, we have 
\begin{equation}\label{LP2}
s\max\{0,u\cdot v\} = h_{[o,su]}(v)\leq h_{K_m}(v) \quad \text{for all }v\in \mathbb{S}^{n-1}
\end{equation}
From (\ref{LP1}) and (\ref{LP2}), it follows that 
\begin{align*}
C'' s \leq \int\displaylimits_{\mathbb{S}^{n-1}}s \max\{0,u\cdot v\}dS_i(K_m,v) &\leq \int\displaylimits_{\mathbb{S}^{n-1}}h(K_m,v)dS_i(K_m,v)\\  & = nV_{i+1}(K_m)\leq n C' 
\end{align*}
Therefore, $s\leq nC'/C''$. Since $0 \in K_m$ for all $m$, we conclude that, for all $m$, $ K_m \subseteq (nC'/C'')\mathbb{B}^n$. By the Blaschke selection theorem {\bf\cite[\textnormal{Theorem 1.8.7}]{schneider93}}, every subsequence of $(K_m)_{m=1}^\infty$ has a subsequence that converges in the Hausdorff metric. The weak convergence of the area measures and Aleksandrov--Fenchel--Jessen theorem \textbf{\cite[\textnormal{Theorem 5.2.3}]{schneider93}}  imply that such subsequences converge to $K$ with respect to the Hausdorff metric and so $d_H(K_m,K)\to 0$ as $m\to \infty$. The proof of the reverse implication is immediate. 

\hfill $\blacksquare$

The following lemma is an extension of Lemma 3 point (3) in \textbf{\cite{FNRZ2011}} to a larger class of multiplier operators.

\begin{lem}\label{smoothing} Let $g$ be a zonal function on $\mathbb{S}^{n-1}$ such that $a_k^n[g] = O(k^{-\beta})$. If $f\in \mathcal{U_\alpha}$, then $\mathrm{T}_g f \in \mathcal{U}_{\alpha + \beta}$ and there is a constant $C_{g}$ such that $\|\mathrm{T}_g f\|_{\mathcal{U}_{\alpha + \beta}}\leq C_{g}\|f\|_{\mathcal{U}_\alpha} $.
\end{lem}
\noindent {\it Proof.} By Parseval's identity  and (\ref{Mn}), we have that 
\begin{align*}
\lVert \mathrm{T}_g f - \mathrm{T}_g\mathrm{M}_j f \rVert^2_{2}  &= \sum_{k = j}^\infty a_k^n[g]^2(1-\Theta(k/j))^2 \rVert \pi_k f\lVert^2_{2}\\ 
&\leq C j^{-2 \beta} \sum_{k = j}^\infty (1-\Theta(k/j))^2 \rVert \pi_k f\lVert^2_{2}\\
& \leq C j^{-2 \beta} \lVert f - \mathrm{M}_j f \rVert^2_{2}\\
& \leq C j^{-2 (\beta +\alpha)}\|f\|_\alpha^2,
\end{align*}
where, in the second line, we use the fact that $\Theta(k/j)\in[0,1]$ and, in the last inequality, we use Lemma~\ref{Uclassesmult0}.
\hfill $\blacksquare$

\vspace{0.3cm} 

The next lemma was established recently in \textbf{\cite{OrtegaSchuster2021}}. 

\begin{lem} \label{Cphi} Suppose that $g \in C^2([-1,1])$. Define $h \in C^2(\mathbb{R}^n\setminus \{0\})\cap  C(\mathbb{R}^n)$ by
\[h(x) = \left \{ \begin{array}{ll} |x|\, g \hspace{-0.1cm} \left( \displaystyle{x \cdot \bar{e}/|x| } \right ) & \mbox{ for all } x \neq o, \\ 0 & \mbox{ for } x = o.    \end{array}    \right .  \]
Then
\begin{equation} \label{hessian1717}
D^2h (u) =  \left(g(u \cdot \bar{e})-(u\cdot \bar{e})\phi'(u\cdot \bar{e})\right)\mathrm{p}_{u^\bot} + g''(u\cdot \bar{e}) (\mathrm{p}_{u^\bot}\bar{e} \otimes \mathrm{p}_{u^\bot}\bar{e}),
\end{equation}
where $\mathrm{p}_{u^\bot}= \mathrm{Id} - u\otimes u$ denotes the orthogonal projection onto the hyperplane $u^\bot$. \end{lem}

The following theorem is one of the main results in \textbf{\cite{OrtegaSchuster2021}}. It will be a key ingredient in the proof of our main theorems. It provides a spectral gap for the multipliers of convolution transforms defined in terms of symmetric convex bodies of revolution.

\begin{theorem} \label{spectralgap} Suppose that $L \in \mathcal{K}^n$ is origin-symmetric and $\mathrm{SO}(n - 1)$ invariant. Then
\begin{equation*}
|a_k^n[L]| < \frac{a_0^n[L]}{(k-1)(n+k-1)}
\end{equation*}
for every $k > 2$ and
\begin{equation*}
|a_2^n[L]| \leq \frac{a_0^n[L]}{n+1},
\end{equation*}
where this inequality is also strict if $L$ is of class $C^2_+$.
\end{theorem}

We will also need the asymptotic behavior of the multipliers of $\mathrm{T}_g$ for a smooth function $g$. This is the content of our next lemma.

\begin{lem}\label{multide} Let $m\geq 2$ be an even integer and let $g$ be a zonal function on $\mathbb{S}^{n-1}$ of class $C^m$. Then
\begin{equation}\label{asysmooth}
 a_k^n[g] = \mathrm{O}\!\left( k^{-m - \frac{n-2}{2}} \right) 
\end{equation}
as $k\to \infty$.
\end{lem}

\noindent {\it Proof.} In order to show (\ref{asysmooth}), note that, by (\ref{multconvpkn}), (\ref{deltasmult}), and the fact that the spherical Laplacian $\Delta_{\mathbb{S}}$ is self-adjoint,
\begin{equation*}
a_k^n[g]  = \frac{(-1)^\frac{m}{2}}{k^\frac{m}{2}(k+n-2)^\frac{m}{2}} \int\displaylimits_{\mathbb{S}^{n-1}}\!\! \Delta_{\mathbb{S}}^\frac{m}{2}g(u)P_k^n(\bar{e} \cdot u)\,du.
\end{equation*}
Hence, by the Cauchy--Schwarz inequality,
 \begin{equation*}
|a_k^n[g]|  \leq \frac{1}{k^m}\sqrt{\frac{\omega_n}{N(n,k)}}{\left\lVert\Delta^\frac{m}{2}_{\mathbb{S}}g\right\rVert}_{L^2} < \infty.		
\end{equation*}
Consequently, by (\ref{nnk}), we obtain the desired asymptotic estimate
\begin{equation*}
a_k^n[g] = \mathrm{O} \left( k^{-m - \frac{n-2}{2}} \right)  \quad \mbox{as } \,\,k \to \infty.
\end{equation*}

\vspace{-0.5cm}

\hfill $\blacksquare$

\vspace{0.3cm} 
 
The following lemma shows that the multipliers of $\mathrm{T}_g$ for a given zonal function $g$ on $\mathbb{S}^{n-1}$ can be used to compute those of its derivatives if the dimension is high enough.

\begin{lem}\label{multide} Let $j$ be positive integers such that  $n \geq 2(j+1)$ and $g$ be a zonal function on $\mathbb{S}^{n-1}$ of class $C^j$. Then, 
\begin{equation}
a_k^n\left[g^{(\,j)}\right] = (2\pi)^j a_{k+j}^{n-2j}[g]
\end{equation}
for all $k\geq 0$.
\end{lem}
\noindent {\it Proof.} First, we prove the case $j = 1$. Using formula (\ref{multleg}) and integration by parts yields  
\begin{align*}
\frac{a_{k+1}^{n-2}[g]}{\omega_{n-3}}  	&= \int\displaylimits_{\!-1}^1g(t)P^{n-2}_{k+1}(t)(1-t^2)^\frac{n-5}{2}dt\\
					%&= (-1)^{k+1} \frac{\Gamma\left(\frac{n - 3}{2}\right)}{2^{k+1} \Gamma\left( \frac{n - 1}{2} + k \right)}\int\displaylimits_{\!-1}^1g(t)\left(\frac{d}{dt}\right)^{k+1}(1-t^2)^{\frac{n-3}{2} + k}dt\\
					%&= (-1)^{k} \frac{\Gamma\left(\frac{n - 3}{2}\right)}{2^{k+1} \Gamma\left( \frac{n - 1}{2} + k\right)} \int\displaylimits_{\!-1}^1g'(t)\left(\frac{d}{dt}\right)^{k}(1-t^2)^{\frac{n-3}{2} + k}dt\\
					&=  \frac{1}{n-3}\int\displaylimits_{\!-1}^1 g'(t)P_{k}^n(t)(1-t^2)^\frac{n-3}{2}dt =  \frac{1}{n-3}\frac{a_k^n[g']}{\omega_{n-1}}.
\end{align*}
Since $ \omega_{n-3}  =  \frac{n-3}{2\pi} \, \omega_{n-1}$, we obtain, 
\[  a_k^n[g'] = 2\pi \, a_{k+1}^{n-2}[g].\]
The general case follows by a simple inductive argument. 

\hfill $\blacksquare$

\vspace{0.3cm}

Now we turn to studying the behavior of integral operators when applied to functions of the $\mathcal{U}_\alpha$ classes. Later, we will see that the hypothesis of the following lemma is satisfied by a variety of generating functions $g$. However, we will leave this as one of the final steps of the proof of Theorem~\ref{mainthm}.   

\begin{lem}\label{smoothing} Let $n \geq 3$, $\alpha\geq 0$, and $g \in C^2([-1,1])$. Suppose that there exists $\beta>0$, and $C_{g, \alpha,\beta}>0$ such that $$\|\nabla^2_{ij} \mathrm{T}_gf\|_{\mathcal{U}_{\alpha +\beta}} \leq C_{g,\alpha,\beta}\|f\|_{\mathcal{U}_\alpha}, \quad i,j\in\{1,\dots, n\}.$$
If $f_1,\dots f_n \in \mathcal{U}_\alpha$, then $\mathrm{D}(D^2\mathrm{T}_g f_1,\dots,D^2\mathrm{T}_g f_{n-1})\in \mathcal{U}_{\alpha+\beta }$ and there exists a constant $C'_{g,\alpha,\beta}>0$ such that 
 \begin{equation}
  \rVert \mathrm{D}(D^2\mathrm{T}_g f_1,\dots,D^2\mathrm{T}_g f_{n-1}) \lVert_{U^{\alpha + \beta }} \leq C'_{g,\alpha,\beta} \rVert f \lVert_{U^{\alpha}}\cdots \rVert f_{n-1} \lVert_{U^{\alpha}}
 \end{equation}
\end{lem}

\noindent {\it Proof.} First we note that for all $u\in \mathbb{S}^{n-1}$,
\begin{equation}\label{dimnsiontrick}
\mathrm{D}(D^2\mathrm{T}_g f_1|_u,\dots,D^2\mathrm{T}_g f_{n-1}|_u) = n \mathrm{D}\left(u \otimes u, \nabla^2 \mathrm{T}_g f_1|_u, \cdots, \nabla^2 \mathrm{T}_g f_{n-1}|_u\right).
\end{equation}
The latter equality follows easily by expanding the mixed discriminant on the right hand side according to the first entry. Hence, 
\begin{equation*}
\mathrm{D}(D^2\mathrm{T}_g f_1,\dots,D^2\mathrm{T}_g f_{n-1}) = \frac{1}{(n-1)!}\sum_{\delta,\sigma\in S_n} \!(-1)^{\mathrm{sgn}(\delta) + \mathrm{sgn}(\sigma)} u_{\delta_1} u_{\sigma_1} \prod_{k=1}^{n-1}\nabla_{\delta_2\sigma_2}(\mathrm{T}_{g} f_k).
\end{equation*}
By assumption $\|\nabla_{ij} \mathrm{T}_gf\|_{\mathcal{U}_{\alpha + \beta}} \leq C \|f\|_{\mathcal{U}_\alpha}$ for all $i,j\in \{1,\dots,n-1\}$. Thus, by the triangle inequality and Lemma~\ref{Uclassesmult}, we obtain 
\begin{align*}
\|\mathrm{D}(D^2\mathrm{T}_g f_1,\dots,D^2\mathrm{T}_g f_{n-1})\|_{\mathcal{U}_{\alpha + \beta}}  \!\!& \leq C'\!\!\!\sum_{\delta,\sigma\in S_n}\left\|(-1)^{\mathrm{sgn}(\delta) + \mathrm{sgn}(\sigma)}u_{\delta_1}u_{\sigma_1} \!\!\prod_{k=1}^{n-1}\nabla_{\delta_2\sigma_2}(\mathrm{T}_{g} f_k)\right\|_{\mathcal{U}_{\alpha + \beta}} \\
& \leq C'' \max_{i,j}\prod_{k = 1}^{n-1}\|\nabla_{ij}\mathrm{T}_{g} f_k\|_{\mathcal{U}_{\alpha + \beta}} \\
& \leq C''' \prod_{k = ^1}^{n-1}\|f_k\|_{\mathcal{U}_\alpha}.
\end{align*}
where $C',C'',C'''>0$ are constant depending only on $g$, $\alpha,\beta$.

\hfill $\blacksquare$

\vspace{0.3cm}

Next we prove some continuity properties of the multilinear map defined by $(f_1,\dots,f_{n-1})\to \mathrm{D}(D^2\mathrm{T}_L f_1,\dots,D^2\mathrm{T}_L f_{n-1})$ for $f_1,\dots,f_{n-1}\in L^2(\mathbb{S}^{n-1})$.

\begin{lem}\label{multilinearbaound} Let $n \geq 3$, and $L$ be a convex body of class  $C_+^2([-1,1])$. Then, 
\begin{itemize}
\item[$(a)$] $\|\mathrm{D}(D^2\mathrm{T}_L f_1,\dots,D^2\mathrm{T}_L f_{n-1})\|_2 \leq \frac{a_0^n[L]}{n-1}\,\|f_1\|_2\prod_{k = 2}^{n-1}\|f_k\|_\infty$,
\item[$(b)$] $\|\mathrm{D}(D^2\mathrm{T}_L f_1,\dots,D^2\mathrm{T}_L f_{n-1})\|_\infty \leq (a^n_0[L])^{n-1} \prod_{k = 1}^{n-1}\|f_k\|_\infty,$
\end{itemize}
for all $f_1\dots,f_{n-1}\in C(\mathbb{S}^{n-1})$. 
\end{lem}
\noindent {\it Proof.} We normalize $T_L$ such that $a^n_0[L] = 1$ (note that this is always possible since $a^n_0[L]>0$). For each $v \in \mathbb{S}^{n-1}$, let $h_{L(v)}$ be the rotated copy of $L$ with axes of revolution $v\in\mathbb{S}^{n-1}$. Since $h_L\in C^2(\mathbb{S}^{n-1})$, it can be easily verified that
\[D^2 \mathrm{T}_L f (u) = \int\displaylimits_{\mathbb{S}^{n-1}}D^2h_{L(v)}(u) f(v)dv,  \qquad u \in \mathbb{S}^{n-1}.\] 
Hence, using the linearity and continuity of the mixed discriminant in each coordinate, we obtain 
\begin{multline}\label{intformula}
\mathrm{D}(D^2\mathrm{T}_L f_1,\dots,D^2\mathrm{T}_L f_{n-1}) = \!\!\!\!\!\!\int\displaylimits_{(\mathbb{S}^{n-1})^{n-1}}\!\!\!\!\!\! \mathrm{D}(D^2 h_{L(v_1)},\dots,D^2 h_{L(v_{n-1})})\prod_{k=1}^{n-1}f_k(v_k)dv,
\end{multline}
where $dv = dv_1 \cdots dv_{n-1}$ is the product measure. Hence, since the mixed discriminant $\mathrm{D}(D^2 h_{L(v_1)},\dots,D^2 h_{L(v_{n-1})})$ is positive, 
\begin{align*}
|\mathrm{D}(D^2\mathrm{T}_L f_1,\dots,D^2\mathrm{T}_L f_{n-1})| 
 &\leq \int\displaylimits_{(\mathbb{S}^{n-1})^{n-1}}\mathrm{D}(D^2 h_{L(v_1)},\dots,D^2 h_{L(v_{n-1})})\prod_{k=1}^{n-1}|f_k(v_k)|dv\\
 &\leq \prod_{k=2}^{n-1}\|f_{k}\|_\infty \!\!\!\!\!\!\!\!\int\displaylimits_{(\mathbb{S}^{n-1})^{n-1}}\!\!\!\!\!\!\mathrm{D}(D^2 h_{L(v_1)},\dots,D^2 h_{L(v_{n-1})})|f_1(v_1)|dv.
\end{align*}
Since we normalized $\mathrm{T}_L$ such that $\mathrm{T}_L h_{\mathbb{B}^{n}} = h_{\mathbb{B}^{n}}$ and since $D^2 h_{\mathbb{B}^{n}} = \mathrm{Id}$,  (\ref{intformula}) yields 
\begin{equation}\label{int1}
\int\displaylimits_{(\mathbb{S}^{n-1})^{n-1}} \mathrm{D}(D^2 h_{L(v_1)},\dots,D^2 h_{L(v_{n-1})})|f_1(v_1)|dv = \mathrm{D}(D^2\mathrm{T}_L |f_1|,\mathrm{Id}[n-2] )
\end{equation}
On the other hand, (\ref{int1}) together with (\ref{basicMD}) yields 
\begin{equation*}
\int\displaylimits_{(\mathbb{S}^{n-1})^{n-1}} \mathrm{D}(D^2 h_{L(v_1)},\dots,D^2 h_{L(v_{n-1})})|f_1(v_1)|dv = \frac{1}{n-1}\mathrm{tr}(D^2 T_L |f_1|) = \Box_n \mathrm{T}_L |f_1|,
\end{equation*}
since for all $f\in C^2(\mathbb{S}^{n-1})$ we have that $D^2 f = \nabla^2_\mathbb{S}f + f I$ and thus $\frac{1}{n-1}\mathrm{tr}(D^2f) = \frac{1}{n-1}\Delta_\mathbb{S}f + f$. Hence, we have that  
\begin{equation}\label{fast}
|\mathrm{D}(D^2 T_{L}f_1,\dots,D^2 T_{L}f_{n-1})| \leq |\Box_n \mathrm{T}_L |f_1||\prod_{k = 2}^{n-1}\|f_k\|_\infty
\end{equation}
Therefore, taking the $L^2$ norm of both sides of inequality (\ref{fast}) we obtain 
\begin{align*}
\|\mathrm{D}(D^2\mathrm{T}_L f_1,\dots,D^2\mathrm{T}_L f_{n-1})\|_2 &\leq \|\Box_n \mathrm{T}_L |f_1|\|_2\prod_{k = 2}^{n-1}\|f_k\|_\infty\\
&\leq  \sup_{k\geq 0}{a_k^n[\Box_n T_L]}  \|f_1\|_2\prod_{k = 2}^{n-1}\|f_k\|_\infty.
\end{align*}
Finally, by Theorem~\ref{spectralgap} we have that 
$$\sup_{k\geq 0}{|a_k^n[\Box_n T_L]|} = \sup_{k\geq 0}{\frac{(k-1)(k+n-1)|a_k^n[L]|}{n-1}} \leq \frac{a_0^n[L]}{n-1}.$$ 
This concludes the proof of $(a)$. Part $(b)$ follows immediately from (\ref{intformula}).

\hfill $\blacksquare$

\vspace{0.4cm} 

\reseteqn \alpheqn \setcounter{theorem}{0}

\vspace{0.6cm}

\centerline{\large{\bf{ \setcounter{abschnitt}{6}
\arabic{abschnitt}. Proof of the main result}}}

\reseteqn \alpheqn \setcounter{theorem}{0}

\vspace{0.6cm}

In this section we prove the main results of this article. We start by showing the global convergence of the iterations of smooth Minkowski valuations of degree $1$. It serves as a guideline to understand the higher degree cases.

\vspace{0.6cm}

\noindent {\it Proof of Theorem~\ref{mainthm2}.} Since $\Phi_1$ is non-trivial, we have that $a_0^n[g] >0$. We may normalize $\Phi_1$ such that  $a_0^n[g] = 1$. Let $K \in \mathcal{K}^n$ and let $\gamma = 1/\pi_0 h_K$. Parseval's identity implies that  
\begin{equation}\label{eqeq1}
\|h_{\Phi_1 \gamma K}- h_{\mathbb{B}^n}\|_2^2 = \sum_{k = 2}^\infty a_k^n[\Box_n g]^2\|\pi_k h(\gamma K,\cdot)\|_2^2 .
\end{equation} 
On the other hand, note that by (\ref{boxhks1}), for any convex body $L$,
\begin{equation}
h_{\Phi_1 L} = S_1(L,\cdot)\ast g = \Box_n h_L \ast g = h_L \ast \Box_n g, 
\end{equation} 
hence, $\Box_ng$ is the generating function of $\Phi_1$ when considered as Minkowski endomorphism. Since $\Phi_1$ is assumed to be monotone we have that $ f = \Box_ng\geq 0$ (see Theorem 1.3 part (iii) in \textbf{\cite{kiderlen05}}). Moreover, since for all $k\geq 1$ and $n\geq 3$ the Legendre polynomial $P_k^n$ satisfy $|P_n^k(t)|<1$ for  all $t\in (-1,1)$ (see Lemma~2.2 in \textbf{\cite{kiderlen05}}),
\begin{equation} \label{Grogu}
|a_k^n[f]| \leq \omega_{n-1}\displaystyle{\int_{-1}^1f(t)|P_k^n(t)|(1-t^2)^\frac{n-3}{2}dt < \omega_{n-1}\displaystyle\int_{-1}^1f(t)(1-t^2)^\frac{n-3}{2}dt = a_0^n[f]}
\end{equation}
holds for all $k \geq 2$ since $\Phi_1$ is non-trivial and monotone. Moreover, by (\ref{multide}), 
$$
a_k^n[f] = a_k^n[\Box_ng] = \frac{(k-1)(k+n-1)}{n-1}a_k^n[g] = \mathrm{O}(k^{-\frac{n-2}{2}})\to 0 \quad \text{as }k\to \infty.
$$
Thus, we obtain 
$$\lambda_g = \sup_{k\geq 2}|a_k^n[f]|< a_0^n[f] = a_0^n[\Box_n]a_0^n[g] = 1.$$
Hence, inequality (\ref{eqeq1}) yields   
\begin{equation}\label{contraction1}
\|h_{\Phi_1 \gamma K}- h_{\mathbb{B}^n}\|_2 \leq \lambda_g \|h_{\gamma K}- h_{\mathbb{B}^n}\|_2
\end{equation}
Therefore, by iterating (\ref{contraction1}) and taking limit as $m\to \infty$ we deduce that 
\begin{equation}
\|h(\Phi_1^m \gamma K,\cdot)- h(\mathbb{B}^n,\cdot)\|_2 \leq \lambda_g^m  \|h(\gamma K,\cdot)- h(\mathbb{B}^n,\cdot)\|_2 \to 0 \quad \text{as } m \to \infty,
\end{equation} 
since $\lambda_g \in (0,1)$. In other words, 
$$d_2(\gamma \Phi^m_1 K,\mathbb{B}^n)\to 0 \quad \text{as } m \to \infty.$$
Since $d_2$ and $d_H$ generate the same topology on $\mathcal{K}^n$ (see Theorem 3 in \textbf{\cite{Vitale1985}}), we conclude that $d_H(\gamma \Phi^m_1 K,\mathbb{B}^n)\to 0$ as $m \to \infty.$
\hfill $\blacksquare$

\vspace{0.3cm}

Now we prove a general theorem from which Theorem \ref{mainthm} will be derived easily.  From now on, we normalize $T_L$ such that its first multiplier is always $1$ (we can always do this since $a_k^n[L]>0$).

\begin{theorem}\label{generalthm}
Let $2\leq i \leq n-1$ and $\Phi_i: \mathcal{K}^n \to \mathcal{K}^n$ be a continuous translation invariant Minkowski valuation of degree $i$ which commutes with $SO(n)$ and suppose that its generating function is the support function of a convex body of revolution $L$ satisfying 
\begin{itemize}
\item[$(a)$] $\mathrm{T}_{L} f \in C^2(\mathbb{S}^{n-1})$ whenever $f\in C(\mathbb{S}^{n-1})$, 
\item[$(b)$] for all $k\geq 2$, the multipliers of $T_L$ satisfy $(k-1)(k+1-2)|a_k^n[L]| \leq (n-1)i^{-1} $, and
\item[$(c)$] there exists $\beta>0$ such that for all $\alpha\geq 0$ and $f \in \mathcal{U}_\alpha$, there exists a constant $C_L = C_{L,\beta,\alpha} >0$ such that  
$$\|\nabla^2_{ij} \mathrm{T}_Lf\|_{\mathcal{U}_{\alpha +\beta}} \leq C_L\|f\|_{\mathcal{U}_\alpha}, \quad i,j\in\{1,\dots, n\}.$$
\end{itemize}
Then, there exists $ \varepsilon >0$ such that if $K \in \mathcal{K}^n$ has continuous surface area measure of order $i$ with density $s_{i}(K,\cdot)$ such that $\|1-s_{i}(\gamma K,\cdot)\|_\infty < \varepsilon$ for some $\gamma>0$, then there exists a sequence of positive numbers $(\gamma_m)_{m = 0}^\infty$ such that  
\begin{equation}
d_{TV}(S_i(\gamma_m\Phi_i^m K,\cdot\,),\sigma_{n-1})\to 0 \quad \text{as }\,\, m \to \infty.
\end{equation}
\end{theorem}

Before we move to the proof of Theorem \ref{generalthm}, we will need to establish the following crucial iteration lemma.  

\begin{lem}\label{iteration step}
Let $\alpha,C_1,C_2>0$ be the constants in Lemma~\ref{fromL2toLinf} and let $L$ be a body of revolution satisfying conditions $(a)$, $(b)$, and $(c)$. There exists $\varepsilon_{L}>0$ and $\lambda_L \in (0,1)$ with the following property: for every $\varepsilon \in (0,\varepsilon_L)$ and every function $f$ such that $\pi_0f = 1$, $\lVert f- 1 \rVert_{L^2}\leq \varepsilon$ and $\lVert f-1 \rVert_{U_\alpha}\leq C_1$, there exists a positive number $\gamma$ such that 
$$\tilde{f} = \gamma \mathrm{D}(D^2 \mathrm{T}_L f [i],\mathrm{Id}[n-1-i])$$  
satisfies $\pi_0 \tilde
{f} = 1$, $\lVert \tilde{f}-1 \rVert_{L^2}\leq \lambda_L \,\varepsilon$ and $\lVert \tilde{f}-1 \rVert_{U_\alpha}\leq C_1$.
\end{lem}

\noindent {\it Proof of Lemma \ref{iteration step}.} Fix $\varepsilon_L>0$. We will determine the conditions needed for $ \varepsilon_L>0$ to satisfy Lemma \ref{iteration step}. Let $\varepsilon \in (0,\varepsilon_L)$ and let $f$ be a function satisfying  $\pi_0f = 1$, $\lVert f- 1 \rVert_{L^2}\leq \varepsilon$ and $\lVert f-1 \rVert_{U_\alpha}\leq C_1$. First note that by the multilinearity of the mixed discriminant and (\ref{basicMD}), we have that 
\begin{align*}
\mathrm{D}(D^2 \mathrm{T}_L f [i],\mathrm{Id}[n-1-i]) &=  \mathrm{D}((D^2 \mathrm{T}_L (f-1) + \mathrm{Id})[i] ,\mathrm{Id}[n-1-i])\\
&= 1 + i \,\,\mathrm{D}(D^2\mathrm{T}_L (f-1),\mathrm{Id}[n-1]) + r \\
&= 1 + i \,\,\mathrm{tr}(D^2\mathrm{T}_L (f-1)) + r\\
&= 1 + i \,\,\Box_n\mathrm{T}_L (f-1) + r,
\end{align*}
where $$r = \sum_{k = 2}^i \frac{i!}{k!(i-k)!}\mathrm{D}(D^2\mathrm{T}_L (f-1)[k],\mathrm{Id}[n-1-k]).$$ 
Notice that by Lemma~\ref{fromL2toLinf}, we have that 
\begin{equation}\label{norminf}
\|f-1\|_{\infty}< C_2 \varepsilon^\frac{4}{n+3}.
\end{equation}
Hence, by Lemma~\ref{multilinearbaound} part $(a)$ and (\ref{norminf}), it follows that  
\begin{align}\label{r2}
\begin{split}
\|r - \pi_{0}r\|_{L^2}  \leq \|r\|_{L^2} & \leq \sum_{k = 2}^i \frac{i!}{k!(i-k)!}\|\mathrm{D}(D^2\mathrm{T}_L (f-1)[k],\mathrm{Id}[n-1-k])\|_{L^2}\\
&\leq   \sum_{k = 2}^i \frac{i!}{k!(i-k)!}\|f-1\|^{k-1}_{\infty}\|f-1\|_{L^2}\\
&\leq   \sum_{k = 2}^i \frac{i!}{k!(i-k)!}C_2^{k-1} \varepsilon^{\frac{4(k-1)}{n+3}+1}\\
&\leq   2^iC_2^{i-1}\varepsilon^\frac{n+7}{n+3},\\
\end{split}
\end{align}
provided that $\varepsilon_L < 1$. In the same fashion, by applying Lemma~\ref{multilinearbaound} part $(b)$, we obtain 
\begin{equation}\label{rinf}
\|r-\pi_0 r\|_\infty \leq 2\|r\|_\infty \leq  2 \sum_{k = 2}^i \frac{i!}{k!(i-k)!}\|f-1\|^{k}_{\infty} \leq 2^{i+1}C_2^i \varepsilon^\frac{8}{n+3}.
\end{equation}
Let $\gamma = \frac{1}{1+ \pi_0r} $, and note that 
\begin{align*}
\tilde{f} = \gamma \mathrm{D}(D^2 \mathrm{T}_L f [i],\mathrm{Id}[n-1-i]) = 1 + \gamma( i\, \Box \mathrm{T}_L(f-1) + r - \pi_0 r).
\end{align*}
By Parseval's identity,
\begin{align*}
\|\Box_n\mathrm{T}_L(f-1)\|_2^2 &= \sum_{k = 1}^\infty \frac{(k-1)^2(k+1-2)^2}{(n-1)^2}a_k^n[L]^2\|\pi_kf\|^2_{2}\\
&\leq \Lambda_L^2 \sum_{k = 1}^\infty \|\pi_k f\|^2_{2} = \Lambda_L^2 \|f-1\|_{2}^2.
\end{align*}
where $$\Lambda_L = \sup_{k\geq 0}\frac{(k-1)(k+1-2)}{n-1}|a_k^n[L]|.$$ 
We obtain,
\begin{equation}\label{contract}
\|\Box_n\mathrm{T}_L(f-1)\|_2 \leq \Lambda_L \|f-1\|_{2}.
\end{equation}
Note that, by Lemma~\ref{asysmooth} and condition (b) in Theorem~\ref{generalthm}, $\Lambda_L <\frac{1}{i}.$ If $\varepsilon_L>0$ is such that $2^{i+1}C_2^{i}\varepsilon_L^\frac{8}{n+3}\leq \frac{1}{4}$, then 
\begin{equation}
0<\gamma = \frac{1}{1 + \pi_0 r}  = 1 - \pi_0 r + \frac{\pi_0r^2}{1 + \pi_0 r} \leq 1 +2^{i+3}C_2^{i}\varepsilon^\frac{8}{n+3},
\end{equation}
since  $|\pi_0r|\leq \|r\|_\infty\leq 2^{i+1}C_2^{i}\varepsilon^\frac{8}{n+3}\leq \frac{1}{4}$. Now fix a number $\lambda_L$ such that $i \Lambda_L<\lambda_L<1$. Then, using (\ref{r2}) and (\ref{contract}) we obtain
\begin{align*}
\lVert\tilde{f}-1\rVert_{2} &\leq \gamma (i\lVert\Box_n \mathrm{T}_L (f-1)\rVert_{2} + \lVert r-\pi_0r\rVert_{2})\\
&\leq (1 + 2^{i+3}C_2^{i}\varepsilon^\frac{8}{n+3})(i\Lambda_L \varepsilon +  2^iC_2^{i-1}\,\,\varepsilon^\frac{n+7}{n+3})\\
&\leq (i\Lambda_L + 2^{i+3}C_2^{i} \varepsilon^\frac{4}{n+3} +   2^iC_2^{i-1} \varepsilon^\frac{8}{n+3} + 2^{2i+3}C_2^{2i-1} \varepsilon^\frac{12}{n+3})\varepsilon\\
&\leq (i\Lambda_L + 2^{2i+4}C_2^{2i-1} \varepsilon^\frac{4}{n+3})\varepsilon\\
&\leq \lambda_L \,\varepsilon,  
\end{align*}
provided $\varepsilon_L>0$ satisfies $2^{2i+4}C_2^{2i-1} \varepsilon_L^\frac{4}{n+3} \leq \lambda_L - i\Lambda_L$.
On the other hand, by (\ref{norminf}) and (\ref{rinf}),
\begin{align}\label{finf}
\begin{split}
\|\tilde{f}-1\|_\infty &\leq \gamma\, (i \, \|\Box_n\mathrm{T}_L(f-1)\|_\infty + \|r- \pi_0 r\|_{\infty})\\
&\leq \gamma\, (i \, \|f-1\|_\infty + \|r- \pi_0 r\|_{\infty})\\
&\leq( 1 + 2^{i+3}C_2^{i}\varepsilon^\frac{8}{n+3})\, (i \, C_2 \varepsilon^\frac{4}{n+3} + 2^iC_2^{i-1}\varepsilon^\frac{n+7}{n+3})\\
&\leq 2^{2i+5}C_2^{2i-1} \varepsilon^\frac{4}{n+3}
\end{split}
\end{align}
Finally, note that by Lemma~\ref{smoothing},
\begin{equation*}
\|\tilde{f}-1\|_{\mathcal{U}_{\alpha +\beta}} \leq \|\tilde{f}\|_{\mathcal{U}_{\alpha +\beta}} +1 \leq C_L \|f\|_{\mathcal{U}_{\alpha}}^i+1 \leq C_L(\|f-1\|_{\mathcal{U}_{\alpha}}+ 1)^i + 1 \leq 4 C_L C_1^i
\end{equation*}
Fix $\delta>0$ such that $\delta C_L C_1^i\leq \frac{1}{8}C_1$. Then, by (\ref{finf}) and Lemma~\ref{tradeoff}, we obtain 
\begin{align*}
\|\tilde{f}-1\|_{\mathcal{U}_\alpha} &\leq C_\delta\|\tilde{f}-1\|_\infty + \delta \|\tilde{f}-1\|_{\mathcal{U}_{\alpha +\beta}}\\
&\leq C_{\delta}2^{2i+5}C_2^{2i-1} \varepsilon^\frac{4}{n+3} + \delta 4 C_L C_1^i \leq C_1,
\end{align*}
provided that $\varepsilon_L>0$ is such that $C_\delta2^{2i+5}C_2^{2i-1} \varepsilon_L^\frac{4}{n+3}\leq \frac{1}{2}C_1$.
\vspace{-0.5cm}

\hfill $\blacksquare$

\vspace{0.3cm} 

With the iteration lemma established, we may now conclude the proof of Theorem~\ref{generalthm}.

\vspace{0.3cm}

\noindent {\it Proof of Theorem \ref{generalthm}.} Fix $\alpha$ be such as Lemma \ref{fromL2toLinf} is satisfied. Let $K\in \mathcal{K}^n$ be a convex body with absolutely continuous area measure of order $i$ such that $\|s_i(K,\cdot)-1\|_{\infty}<\varepsilon$ for some $\varepsilon$ that we will fix later. Thus, $\|s_i(K,\cdot)\|_{\mathcal{U}_0}<\varepsilon +1$. By Lemma~\ref{smoothing}, 
\begin{equation*}
\|s_i(\Phi_i K,\cdot\,)\|_{\mathcal{U}_\beta} = \|\mathrm{D}(D^2\mathrm{T}_L s_i(K,\cdot\,)[i],\mathrm{Id}[n-1-i])\|_{\mathcal{U}_\beta} \leq C_L \|s_i(K,\cdot\,)\|_{\mathcal{U}_0}^i,
\end{equation*}
and thus, for some constant $C'_L>0$,
\begin{equation}\label{b0}
\|s_i(\Phi_i^m K,\cdot\,)\|_{\mathcal{U}_{m \beta}}\leq C'_L \|s_i(K,\cdot\,)\|_{\mathcal{U}_0}^{i^m} \leq C'_L(1+\varepsilon)^{i^m}
\end{equation}
As we already pointed out in the proof of Lemma~\ref{multilinearbaound}, 
\[D^2 \mathrm{T}_L f (u) = \int\displaylimits_{\mathbb{S}^{n-1}}D^2h_{L(v)}(u) f(v)dv,  \qquad u \in \mathbb{S}^{n-1}.\] 
Using the normalization $T_L h_{\mathbb{B}^n} = h_{\mathbb{B}^n}$, we obtain 
\begin{align*}
s_i(\Phi_i K,\cdot)	&= \mathrm{D}\left(D^2\mathrm{T}_L s_i(K,\cdot\,)[i],\mathrm{Id}[n-1-i]\right)\\
					&= \mathrm{D}\left(D^2\mathrm{T}_L s_i(K,\cdot\,)[i],D^2\mathrm{T}_L 1 [n-1-i]\right)\\
&= \!\!\!\int\displaylimits_{(\mathbb{S}^{n-1})^{n-1}}\mathrm{D}(D^2h_{L(v_1)},\dots,D^2 h_{L(v_{n-1})})s_i(K,v_1)\cdots s_i(K,v_{i})dx_1\!\cdots dv_{n-1}.
\end{align*}
Since $1-\varepsilon < s_i(K,u) < 1 + \varepsilon$ for all $u\in\mathbb{S}^{n-1}$, we have that 
\begin{equation}\label{b1}
(1-\varepsilon)^{i^m}< s_i(\Phi^m_i K,u) < (1+\varepsilon)^{i^m},
\end{equation}
for all $u\in\mathbb{S}^{n-1}$. Thus, $f = s_i(\Phi^m_i K,\cdot\,)/\pi_0 s_i(\Phi^m_i K,\cdot\,)$ is such that $\pi_0 f = 1$ and, by (\ref{b0}) and (\ref{b1}), 
\begin{equation}
\|f-1\|_{\mathcal{U}_{m\beta}} \leq 1 + \frac{\|s_i(\Phi^m_i K,\cdot\,)\|_{\mathcal{U}_{m\beta}}}{\pi_0 s_i(\Phi^m_i K,\cdot\,)} < 1 + \frac{C'_L(1+\varepsilon)^{i^m}}{(1-\varepsilon)^{i^m}} \leq C'_L 4^{i^m},
\end{equation}
provided that $\varepsilon<\frac{1}{2}$. Fix $m$ large enough so that $ m\beta>\alpha>0$. Clearly from (\ref{b1}), there exists a constant $C'>0$ such that $\|f-1\|_\infty < C'\varepsilon$. Fix $\delta>0$ such that $\delta  C'_L 4^{i^m}< \frac{1}{2}C_1$ where $C_1$ is the constant in Lemma~\ref{iteration step}. Hence by Lemma~\ref{tradeoff}, we obtain  
\begin{align*}
\|f-1\|_{\mathcal{U}_{\alpha}} &< C_\delta \|f-1\|_{\infty} + \delta \|f-1\|_{\mathcal{U}_{m\beta}}\\
&< C_\delta C'\varepsilon + \delta  C'_L 4^{i^m} \leq C_1
\end{align*}
provided that $\varepsilon>0$ is such that $C_\delta C'\varepsilon \leq \frac{1}{2}C_1$. Therefore, without loss of generality, we may assume that $K$ satisfies the conditions of Lemma~\ref{iteration step} (if not we replace $K$ by $\Phi_i^mK$ for large enough $m$). Hence, there exists a $\gamma_1$ such that the body $K_1 = \gamma_1\Phi_iK$ satisfies the assumptions of Lemma \ref{iteration step} with $\lambda_L \varepsilon$ instead of $\varepsilon$. Applying Lemma~\ref{iteration step} again with $\lambda_L^2 \varepsilon$, we obtain a positive number $\xi_2$ such that $K_2 = \xi_2 \Phi_i K_1 = \xi_2 \gamma_1^i \Phi_i^2 K = \gamma_2 \Phi_i^2 K$ satisfies the assumptions of Lemma \ref{iteration step} with $\lambda_L^2 \varepsilon$, where $\gamma_2 = \gamma_1 \xi_2$. Thus, continuing this process inductively, we find a sequence $(\gamma_m)_{m=1}^\infty$ of positive numbers such that
$$\|s_i(\gamma_m\Phi^m_i K,\cdot\,)-1\|_{L^2}\leq \lambda_L^m \varepsilon,$$
and $\|s_i(\gamma_m\Phi^m_i K,u)\|_{\mathcal{U}_\alpha}\leq C_1$. Hence, Lemma~\ref{fromL2toLinf} yields 
$$d_{\mathrm{TV}}(S_i(\gamma_m\Phi^m_i K,\cdot\,),\sigma_{n-1}) \leq \|s_i(\gamma_m\Phi^m_i K,\cdot\,)-1\|_{\infty}\leq C \lambda_L^m \varepsilon\to 0 \quad \text{as }\, m\to \infty,$$
which is the desired result.
\hfill $\blacksquare$

\vspace{0.3cm}

Note that the assumption of Theorem \ref{mainthm} on the support function of $K$ is stronger than the assumption of Theorem \ref{generalthm} for the area density. The area density of order $i$ of $K$ is given in terms of a mixed discriminate involving the Hessian of $h_K$ and the identity matrix. Hence, if we assume that $h_K$ is close enough to $h_{\mathbb{B}^{n}}$ in the $C^2$ norm then, by continuity, $s_i(K,\cdot)$ will be close to the surface area measure of the sphere $s_i(\mathbb{B}^n,\cdot) = \sigma_{n-1}$. 

\vspace{0.3cm}  

\noindent {\it Proof of Theorem \ref{mainthm}.} By Theorem~\ref{generalthm}, it remains to check that $L$ satisfies the condition $(a)$, $(b)$ and $(c)$. Since $h_L\in C^\infty(\mathbb{S}^{n-1})$, as in the proof of Lemma~\ref{multilinearbaound}, we have that 
\begin{equation}\label{intfrom1}
\nabla^2 \mathrm{T}_L f (u) = \int\displaylimits_{\mathbb{S}^{n-1}}\nabla^2 h_{L(v)}(u) f(v)dv,  \qquad u \in \mathbb{S}^{n-1}.
\end{equation}
This verifies $(a)$. Condition $(b)$ of Theorem~\ref{generalthm} is the content of Theorem~\ref{spectralgap}. Finally, to verify that $(c)$ is satisfied. Note that by Lemma~\ref{Cphi},
\begin{equation} \label{hessian1717}
\nabla h_{L(u)} (u) =  \left(g(u \cdot v)-(u\cdot v)g'(u\cdot v)\right)\mathrm{p}_{u^\bot} + g''(u\cdot v) (\mathrm{p}_{u^\bot}v \otimes \mathrm{p}_{u^\bot}v),
\end{equation}
where $g \in C^\infty([-1,1])$ is such that $h_{L(v)}(u) = g(v\cdot u)$ for all $u,v \in \mathbb{S}^{n-1}$.Note that $\mathrm{T}_L = \mathrm{T}_g$. Thus, using (\ref{Cphi}) in (\ref{hessian1717}) yields
\begin{align*}
\nabla^2 \mathrm{T}_g f (v) &= \int\displaylimits_{\mathbb{S}^{n-1}}g(u \cdot v)\mathrm{p}_{u^\perp}f(u)du- \int\displaylimits_{\mathbb{S}^{n-1}}g'(u\cdot v)(u\cdot v)\mathrm{p}_{u^\bot} f(u)du\\ 
&\,\,\,\,\,\, + \int\displaylimits_{\mathbb{S}^{n-1}} g''(u\cdot v) (\mathrm{p}_{u^\bot}v \otimes \mathrm{p}_{u^\bot}v)f(u)du\\
& = \mathrm{T}_g[\mathrm{p}_{u^\perp}f(u)](v) -  \mathrm{T}_{g'}[ (v \cdot u) \,\,\mathrm{p}_{u^\perp}f(u)](v) + \mathrm{T}_{g''}[(\mathrm{p}_{u^\bot}v \otimes \mathrm{p}_{u^\bot}v)f(u)].
\end{align*}
Therefore for all $i,j\in\{1,\dots,n-1\}$, we can write 
\begin{multline}\label{hessTgf}
\nabla_{ij}T_g f(v) = \mathrm{T}_g[(\delta_{ij}- u_iu_j)f(u)] - T_{g'}[(v\cdot u) (\delta_{ij}- u_iu_j) f(u)]\\
 + \mathrm{T}_{g''}[(v_i- (v\cdot u)u_i)(v_j - (v\cdot u)u_j)f(u)]
\end{multline}
Since, for all $i,j\in\{1,\dots,n-1\}$, the functions $p_1(u,v) = \delta_{ij}- u_iu_j$, $p_2(u,v)  =(v\cdot u) (\delta_{ij}- u_iu_j)$, $p_3(u,v)  = (v_i- (v\cdot u)u_i)(v_j - (v\cdot u)u_j)$ are all polynomials they belong to $\mathcal{U}_\alpha$ for all $\alpha>0$. Furthermore, since $g$ is a smooth function, Lemma~\ref{asysmooth} implies that its multipliers and all the multipliers of all of its derivatives decay as fast as any polynomial. Hence, for any fixed $\beta>0$ we have that  $a_k^n[g^{(i)}] = \mathrm{O}(k^{-\beta})$ as $k\to \infty$ for $i = 0,1,2$. Thus, Lemma~\ref{smoothing} together with Lemma \ref{Uclassesmult} imply that 
\begin{equation}\label{pk}
\|\mathrm{T}_{g^{(i)}}[p_k f]\|_{\mathcal{U}_{\alpha + \beta}} \leq C \|f\|_{\mathcal{U}_{\alpha}}
\end{equation}
for all $i\in \{0,1,2\}$, where $C$ is a positive constant depending on $\|v_iv_j\|_{\mathcal{U}_\alpha}$, $\|v_iv_jv_k\|_{\mathcal{U}_\alpha}$,  $\|v_iv_jv_kv_l\|_{\mathcal{U}_\alpha}$ for $i,j,k,l \in \{1,\dots n\}$ , and $L$. Hence, by (\ref{pk}) and (\ref{hessTgf}) we obtain 
$$\|\nabla_{ij}T_g f\|_{\mathcal{U}_{\alpha+\beta }}\leq C'\|f\|_{\mathcal{U}_{\alpha}}$$ 
for some positive constant $C'$. This verifies condition $(c)$. \\
Therefore, we may apply Theorem~\ref{generalthm} to find an $ \varepsilon >0$ such that if $K \in \mathcal{K}^n$ has continuous surface area measure of order $i$ with density $s_{i}(K,\cdot)$ such that $$\|1-s_{i}(\gamma K,\cdot)\|_\infty < \varepsilon,$$ for some $\gamma>0$, then there exists a sequence of positive numbers $(\gamma_m)_{m = 0}^\infty$ such that  
\begin{equation*}
d_{TV}(S_i(\gamma_m\Phi_i^m K,\cdot\,),\sigma_{n-1})\to 0 \quad \text{as }\,\, m \to \infty.
\end{equation*}
Since convergence with respect to the total variation metric implies convergence with respect to the $d_{LP}$-metric, and convergence with respect to the $d_{LP}$-metric is equivalent to convergence with respect to the Hausdorff metric (Proposition~\ref{DpH}), we have that   
\begin{equation*}
d_{H}(\gamma_m\Phi_i^m K,\mathbb{B}^n)\to 0 \quad \text{as }\,\, m \to \infty.
\end{equation*}
\hfill $\blacksquare$

\vspace{0.3cm} 

Let us point out that it is possible to weaken the regularity assumptions on the support function of $L$ if one assumes that the dimension is higher. In this case we employ Lemma~\ref{multide} to verify condition $(c)$ in Theorem~\ref{generalthm}.

\vspace{0.3cm} 

\noindent {\it Proof of Theorem \ref{mainthm} with $L$ of class $C^2_+$ and $n\geq 6$.} The proof that conditions $(a)$, $(c)$, and $(d)$ are satisfied here is essentially the same as the proof of Theorem~\ref{mainthm}. However, to obtain condition $(b)$, note that by Lemma~\ref{asysmooth}, 
\begin{equation}
 a_k^n[g] = \mathrm{O}\!\left( k^{- \frac{n + 2}{2}} \right)\quad \text{as }k\to \infty. 
\end{equation}
By Lemma~\ref{multide}, $a_k^n[g'] = 2\pi\ a_{k + 1}^{n-2}[g]$ and $a_k^n[g''] = (2\pi)^2 a_{k + 2}^{n-4}[g]$ for all $k\geq 0$. Hence, by applying Lemma~\ref{asysmooth} again, we obtain  
\begin{equation}
 a_k^n[g'] = \mathrm{O}\!\left( k^{- \frac{n}{2}} \right) \text{ and } a_k^n[g''] = \mathrm{O}\!\left( k^{- \frac{n-2}{2}} \right) \quad \text{as }k\to \infty. 
\end{equation}
Finally, using (\ref{hessTgf}) we obtain (in the same way as in the proof of Theorem~\ref{mainthm}), 
\begin{equation}
\|\nabla_{ij}\mathrm{T}_g f\|_{\mathcal{U}_{\alpha+ \frac{n-2}{2}}}\leq C'\|f\|_{\mathcal{U}_{\alpha}},
\end{equation}
 for all $i,j\in\{1,\dots, n\}$. 
The rest of the proof is now identical to the proof of Theorem~\ref{mainthm}.
\vspace{0.3cm} 
\hfill $\blacksquare$
\vspace{0.3cm}

Finally, we finish this article by showing how Theorem~\ref{mainthm} can be use to derive a local solution to Conjecture 1 via the class reduction technique (Proposition~\ref{prop1inequ}). First we prove Corollary 1. 

\vspace{0.3cm} 

\noindent {\it Proof of Corollary 1.}
Suppose $K$ satisfies the hypothesis of Theorem~\ref{mainthm}. By iterating $\Phi_i$, we obtain 
$$\Phi_i^{2m} K = \alpha_m K,$$
for some sequence of  real numbers $\alpha_m$. Therefore, we can find a sequence $(\gamma_{m})_{m = 1}^\infty$ by Theorem 1 such that 
$$\gamma_{2m}\alpha_{m}K = \gamma_{2m}\Phi_i^{2m} K   \to \mathbb{B}^n \quad \text{as } m\to \infty$$ in Hausdorff metric. By continuity of the mean width we derive $\gamma_{2m}\alpha_{m} \to 2/w(K),$ and therefore,   
$\gamma_{2m}\alpha_{m} K \to \frac{2}{w(K)}K$. Hence, by uniqueness of the limit, $K = \frac{w(K)}{2}\mathbb{B}^n.$ 
\vspace{0.3cm} 
\hfill $\blacksquare$

\vspace{0.3cm} 

\noindent {\it Proof of Corollary 2.} Suppose that $K\in \mathcal{K}^n$  satisfies the hypothesis of Theorem~\ref{mainthm}. Let $\psi:\mathcal{K}^n\to \mathbb{R}$ be the functional defined by 
$$\psi(K) = \frac{V_{i+1}(\Phi_iK)}{V_{i+1}(K)^{i}},$$
for all $K\in \mathcal{K}^n$. Since $V_{i+1}$ is homogeneous,
 $$\psi(\gamma_m \Phi_i^m K ) = \psi(\Phi_i^m K)$$ for all $m$. By Proposition~\ref{prop1inequ}, 
\begin{equation}
\psi(K) \geq \psi(\Phi_i K)\geq \cdots \geq  \psi( \Phi_i^m K) \to  \psi(\mathbb{B}^n) \quad \text{as } m\to \infty
\end{equation}
Therefore, $\psi(K)\geq \psi(\mathbb{B}^n)$ as required. To deal with the equality cases locally,  suppose that $\psi(K) = \psi(\mathbb{B}^n)$. The second part of Proposition~1 implies that $\Phi_i^2 K = \alpha K$. By Corollary~1 we obtain the desire conclusion. 
\vspace{0.3cm} 
\hfill $\blacksquare$

\vspace{1cm}
\noindent {{\bf Acknowledgments} The author was supported by the Austrian Science Fund (FWF), Project number: P31448-N35.

\begin{small}
\noindent Oscar Ortega-Moreno\\
Vienna University of Technology\\
oscarortem@gmail.com
\end{small}

\end{document}